\documentclass[final,leqno,onefignum,onetabnum]{siamltex}
\usepackage{amsmath}
\usepackage{amsfonts}
\usepackage{multirow}
\usepackage{graphicx}
\usepackage{mathtools}
\title{A dynamic domain decomposition\\for a class of second order semi-linear equations}

\author{Simone Cacace and Maurizio Falcone}

\begin{document}
\maketitle
\slugger{sinum}{xxxx}{xx}{x}{x--x}

\begin{abstract}
We propose a parallel algorithm for the numerical solution of a class of second order semi-linear equations coming from stochastic 
optimal control problems, by means of 
a dynamic domain decomposition technique. 
The new method is an extension of the patchy domain decomposition method presented in \cite{CCFP12} for first order Hamilton-Jacobi-Bellman 
equations related to deterministic optimal control problems. The semi-Lagrangian scheme underlying the original method is modified in order 
to deal with (possibly degenerate) diffusion, by approximating the stochastic optimal control problem associated to the equation 
via discrete time Markov chains. We show that under suitable conditions on the discretization parameters and for sufficiently 
small values of the diffusion coefficient, the parallel computation on the proposed dynamic decomposition is faster than 
that on a static decomposition. To this end, we combine the parallelization 
with some well known techniques in the context of fast-marching-like methods for first order Hamilton-Jacobi equations. 
Several numerical tests in dimension two are presented, in order to show the features of the proposed method.   

\end{abstract}

\begin{keywords}\end{keywords}

\begin{AMS}\end{AMS}

\pagestyle{myheadings}
\thispagestyle{plain}

\section{Introduction}\label{intro}
Domain decomposition techniques have become popular in the solution of partial differential equations 
arising in several applied 
contexts, including fluid dynamics and electromagnetism.

Given a decomposition of a domain in subsets of manageable size and 
prescribed suitable transmission conditions at the interfaces or overlapping regions between the sub-domains, 
massive parallel computations can be performed exploiting the computing power of modern clusters of CPUs and/or GPUs. 
This is the most common strategy to attack the so called \emph{curse of dimensionality}, a quite general term to 
denote the difficulties related to the numerical approximation of problems in high dimension, 
both in terms of memory requirements and computational efforts.  
New efficient and accurate parallel algorithms are increasingly required, 
to open the way to the solution of real-life problems with a very huge number of degrees of fredom. 

Despite such techniques have been mainly designed, analyzed and implemented in the framework of variational methods for elliptic 
and parabolic equations, in the last decades an increasing interest has also emerged in fields more related to hyperbolic 
equations, e.g. optimal 
control problems, differential games, front propagation and image processing. In this setting, a leading role is played by 
the theory of viscosity solutions for Hamilton-Jacobi equations, representing a solid theoretical background for a growing 
number of numerical methods for the solution of problems in robotics, aeronautics, electrical and aerospace 
engineering. 

In \cite{CCFP12} the authors proposed the patchy domain decomposition method, 
a parallel algorithm for the solution of deterministic optimal control problems, based on a semi-Lagrangian 
discretization of the corresponding first order Hamilton-Jacobi-Bellman equations. 
The main idea there can be summarized as follows. First the solution of the equation is computed on a grid 
which is coarse with respect to the actual (fine) grid. Then, an approximation of the optimal vector field 
associated to the underlying control problem is synthesized on the fine grid. 
Finally, a decomposition of the fine grid is dynamically build driven by 
the approximate optimal vector field, and the equation is solved in parallel on the corresponding sub-domains. 
This pre-computation has clearly a cost and produces a rather complex subdivision of the 
computational box, but gives the fundamental property that each sub-domain is, up to an error, independent on the others. 
This feature allows for an efficient parallelization, since transmission conditions at the interfaces of the sub-domains can be 
completely avoided. The computational resources can be better employed, so that no processor remains idle or computes 
useless information, thus saving CPU time to reach convergence. Despite the pre-computation step and the non trivial implementation, 
the overall performance of the patchy decomposition overcomes that of a static domain decomposition.

Another technique widely used to reduce the computational efforts when computing the solution to first order Hamilton-Jacobi equations, 
is to exploit the so called \emph{causality property}. This term refers to a peculiarity of hyperbolic equations, namely 
the fact that, starting from the boundary data, information propagates in the domain along characteristics at finite speed. 
Fast-marching methods \cite{T95,S99,SV03} have been developed trying to reproduce this feature at the discrete level. 
The main idea is to process the grid nodes in a suitable order that decouples, at least partially, the nonlinear system corresponding to 
the discretization of the equation. Moreover, the computation is localized, at each iteration, 
on the nodes that bring the relevant information. Due to the causality property, these nodes are very few, 
compared to the size of the whole grid. In this way, the solution can be computed in a cascade fashion, 
accelerating the convergence of the underlying scheme significantly. 
It turns out that each node converges in a predetermined and finite number of iterations, 
only one iteration in the most favorable cases. For this reason fast-marching methods are also termed \emph{single-pass} methods.

Keeping these ideas in mind, in this paper we revisit the patchy domain decomposition in the context of second order semi-linear 
equations. We try to extend the main features of the method to a class of problems where diffusion appears, namely stationary 
nonlinear Hamilton-Jacobi-Bellman equations coming from stochastic optimal control problems. 
A prototype problem is the following:
$$
\vspace{-9pt}
\left\{ \begin{array}{ll}
         -\varepsilon\Delta u(x)+\displaystyle\max_{a\in A}\left\{ -f(x,a)\cdot \nabla u(x) \right\} = l(x) & x\in \Omega\\
         u(x)=g(x) & x\in \partial \Omega
        \end{array}
\right.
$$\\
where $\varepsilon>0$ is the diffusion coefficient, $A$ is a compact set of admissible controls, $f$ is the controlled dynamics 
driving the system (aka the \emph{controlled advection}), $l$ is a source term and $g$ is a boundary datum.
 
The semi-Lagrangian local solver of the original patchy method has to be modified in order to deal with diffusion. 
Indeed, the presence of diffusion invalidates the connection with deterministic optimal control, in the sense that characteristic 
curves associated to the hyperbolic equation are no longer well defined in this more general setting. Nevertheless, a weak notion of characteristics 
can be still provided via stochastic differential equations, interpreting the diffusion as a Wiener process. Then, 
a semi-Lagrangian discretization of the corresponding semi-linear equation can be still performed, approximating the stochastic 
process via discrete time Markov chains. The resulting scheme is able to compute the solution of equations with very degenerate 
diffusions. This is a very delicate problem in the community working on advection-diffusion equations, 
especially in cases of interest where the diffusion is much smaller than the advection (e.g. Navier-Stokes equations in fluid dynamics). 
Indeed, it is well known that a degeneracy in the diffusion translates into a lack of coercivity of the variational form 
associated to the equation. 
In particular, discretizations based on finite elements need to be stabilized with ad-hoc procedures, in order to capture 
the boundary or internal layers that typically arise in such equations.

In this more general setting the application of the patchy domain decomposition is not trivial. 
The main obstacle is that the transmission conditions at the interfaces cannot be avoided, due again to the presence of diffusion. 
A natural question 
arises: is the patchy domain decomposition still favorable than a static domain decomposition? We show that the answer to this 
question depends on the ratio between the controlled advection and the diffusion coefficient, 
a characteristic quantity similar to 
the Reynolds number 
for Navier-Stokes equations. More precisely, we show that if the controlled advection dominates diffusion and the discretization parameters are 
chosen in a suitable way, then a parallel computation on the dynamic decomposition can still be faster than that on a static 
decomposition. To reach this achievement, we need to exploit all the information collected in the pre-computation step  
of the dynamic decomposition. In particular, we try to employ, despite the diffusion, the causality property of hyperbolic equations discussed above.
Further technical modifications should be applied to the local solver, in order to remove 
from the scheme the dependency of each node on itself, an issue that dramatically breaks down the causality.
This results in a remarkable speedup for the computation, even more substantial if combined with 
a parallel algorithm.

The paper is organized as follows: in Section 2 we briefly review, for the readers convenience, the patchy domain decomposition 
method for first order Hamilton-Jacobi-Bellman equations proposed in \cite{CCFP12}. Section 3 is devoted to the extension of 
the semi-Lagrangian local solver to second order semi-linear equations, in a form suitable for our purposes. In Section 4, we 
show how to adapt the patchy decomposition method to the more general setting of equations with diffusion, in particular we 
establish conditions on the parameters that guarantee the effectiveness of the dynamic domain decomposition for the parallel 
computation. Finally, in Section 5, we present several numerical tests in dimension two, in order to show the 
performance of the proposed method compared to its static version. 

\section{Patchy decomposition for first order HJB equations}\label{PD1}
In this section we review the patchy domain decomposition method, proposed in \cite{CCFP12} 
to solve boundary value problems for first order Hamilton-Jacobi equations of the form
\begin{equation}\label{HJB1}
\left\{ \begin{array}{ll}
         \displaystyle\max_{a\in A}\left\{ -f(x,a)\cdot \nabla u(x) -l(x,a) \right\} = 0 & x\in \Omega\\
         u(x)=g(x) & x\in \partial \Omega
        \end{array}
\right.
\end{equation}
where $\Omega\subset\mathbb{R}^n$ is an open set and $A\subset\mathbb{R}^m$ is a compact set. Moreover, 
$f:\Omega\times A\to \mathbb{R}^n$ is a vector field and $u:\Omega\to\mathbb{R}$, 
$l:\Omega\times A\to \mathbb{R}$, $g:\partial\Omega\to \mathbb{R}$ are scalar functions.

It is well known that equation $(\ref{HJB1})$ can be interpreted, via the celebrated 
\emph{dynamic programming principle},
as the Hamilton-Jacobi-Bellman equation associated to a suitable deterministic optimal control problem. Indeed, let us consider 
the following controlled dynamical system
\begin{equation}\label{ODE}
\left\{ \begin{array}{ll}
         \dot y(t) = f(y(t),\alpha(t))\,, & t>0\\
         y(0)  = x & 
        \end{array}
\right. 
\end{equation}
where $y(\cdot)$ is the state of the system and $\alpha(\cdot)$ is a generic function belonging to the following set 
of admissible controls: 
$$\mathcal{A}=\{\alpha :[0,+\infty)\to A,\mbox{ measurable}\}\,.$$
We denote by $y(t;x,\alpha(\cdot))$ the solution to $(\ref{ODE})$ starting from $x\in\Omega$ using the control $\alpha(\cdot)$,  
and we define the \emph{first time arrival} to the boundary $\partial \Omega$ as
$$
\tau(x,\alpha(\cdot))=\inf\left\{t\geq 0\,:\, y(t;x,\alpha(\cdot))\in\partial\Omega\right\}\,. 
$$
Finally, we consider the following functional: 
$$
J\big(x,\alpha(\cdot)\big)=\int_0^{\tau(x,\alpha(\cdot))} l\big(y(s;x,\alpha(\cdot)),\alpha(s)\big)\,ds + g\big(y(\tau(x,\alpha(\cdot)))\big)\,.
$$
In this setting, the \emph{minimum time} problem with \emph{running cost} $l$ and \emph{exit cost} $g$ consists in finding, 
for each $x\in\Omega$, an optimal control $\alpha^*(\cdot)\in\mathcal{A}$ that minimizes $J$ among all the admissible controls. 
Under suitable assumptions on the data, it can be proved that the \emph{value function} of the problem, i.e.
$$
u(x)=\inf_{\alpha(\cdot)\in\mathcal{A}} J\big(x,\alpha(\cdot)\big) \qquad x\in\Omega
$$
is the unique viscosity solution to $(\ref{HJB1})$. We refer the reader to \cite{BCD97} for the details.

The main advantage of this approach is that, once the value function $u$ is computed, we can quite easily synthesize an 
optimal control for the minimum time problem, by taking
\begin{equation}\label{synthesis}
a^*(x)=\displaystyle\arg\min_{a\in A}\left\{ f(x,a)\cdot \nabla u(x) +l(x,a) \right\}\qquad x\in\Omega\,.
\end{equation}
Note that $a^*$ is in \emph{feedback form}, namely it depends only on the the state an not on the time. 
It follows that, for each $x\in\Omega$, we can compute the \emph{optimal trajectory} $y^*(\cdot)$ starting at $x$ 
(i.e. a characteristic curve of the hyperbolic equation $(\ref{HJB1})$) by 
simply plugging $a^*$ in the dynamical system $(\ref{ODE})$:
$$
\left\{ \begin{array}{ll}
         \dot y^*(t) = f\left(y^*(t),a^*(y^*(t))\right)\,, & t>0\\
         y^*(0)  = x & 
        \end{array}
\right. 
$$
This is not the case for other types of techniques, based on the \emph{Pontryagin maximum principle}, which provide 
only necessary conditions for optimality and sub-optimal \emph{open-loop} controls (i.e. controls that depend on time). 

Unfortunately, from the numerical point of view, the Hamilton-Jacobi approach implies severe computational efforts, since it requires the computation 
of the value function $u$ on the whole state space $\Omega$. 
Considering that industrial applications demand the solution of optimal control problems at least 
in dimension six (as for the most simple second order dynamics in $\mathbb{R}^3$ with controlled acceleration), this approach still suffers the 
\emph{curse of dimensionality} mentioned in the Introduction. For this reason, the development of efficient parallel algorithms 
for Hamilton-Jacobi equations is nowadays an active field of research. 
The patchy domain decomposition method places itself exactly in this context.

The semi-Lagrangian discretization of equation $(\ref{HJB1})$ is postponed to the next section, 
in a generalized form that recovers first order equations as a particular cases. 
The interested reader can refer to \cite{CCFP12} for the original version of the scheme.
Here, we prefer to keep the discussion free from technical details, 
and mainly focus on the ideas behind the construction of the parallel algorithm. 

We consider two different discretizations $G$ and $G_c$ of the state space $\Omega$. 
The grid $G$ denotes the actual grid on which we want to solve the problem $(\ref{HJB1})$, also termed the \emph{fine} grid, 
whereas $G_c$ is very \emph{coarse} compared to 
(and possibly contained in) $G$. 

The first step of the method consists in computing a coarse solution $u_c$ on $G_c$, also via parallel computations using a standard (static) domain decomposition technique. 
Due to the low resolution of the grid, this is a very cheap operation, but gives a first rough approximation on $G$, say $\hat u$, 
of the actual solution $u$. It is reconstructed by interpolation of $u_c$ on $G$.

Now we employ $\hat u$ to compute a feedback optimal 
control $\hat a^*$, via an appropriate discrete version of the synthesis procedure $(\ref{synthesis})$, 
presented in Section \ref{PD2}. 
We stress that 
$\hat a^*$ is just a coarse approximation of the actual optimal control $a^*$, but it is defined on the fine grid $G$. 
This is enough to start the construction of the patchy domain decomposition. 

We divide the boundary $\partial \Omega$ in a fixed 
number $N_P$ of disjoint sub-sets, denoted by $\Gamma_p$, with $p=1,...,N_P$. Then, for each $p$, 
we compute the sub-domain $\Omega_p\subset \Omega$ as the \emph{numerical domain of dependence} of $\Gamma_p$ through the 
optimal dynamics $f^*(\cdot)=f(\cdot,\hat a^*(\cdot))$. Let us clarify this point in a continuous setting and postpone 
the actual implementation to Section \ref{PD2}. We denote by $\chi_{\Gamma_p}$ the 
characteristic function of the sub-boundary $\Gamma_p$, and consider 
the following Cauchy-Dirichlet problem for the advection equation:
\begin{equation}\label{advection}
\left\{ \begin{array}{ll}
         \partial_t\phi_p(x,t)-f(x,\hat a^*(x))\cdot \nabla \phi_p(x,t)=0 & (x,t)\in\Omega\times (0,+\infty)\\
         \phi_p(x,0)=0 & x\in\Omega \\
         \phi_p(x,t)=\chi_{\Gamma_p}(x) & (x,t)\in\partial\Omega\times (0,+\infty) 
        \end{array}
\right. 
\end{equation}
It is well known that the boundary datum acts as a source of information, that flows (backward) in $\Omega$ 
along characteristics, according to the drift $f^*$. The limit in time 
$$\phi_p^\infty(x)=\displaystyle\lim_{t\to +\infty}\phi_p(x,t)$$ 
is still a characteristic function, since the hyperbolic equation preserves the properties of 
$\chi_{\Gamma_p}$ (e.g. the maximum and the singularities). 
Then, we define the $p$-th patch of our dynamic decomposition as 
$$
\Omega_p=\left\{ x\in\Omega\,:\, \phi_p^\infty(x)=1\right\}\,.
$$

\begin{figure}[htp]
\begin{center}
\begin{tabular}{cc}
 \includegraphics[width=0.45\textwidth]{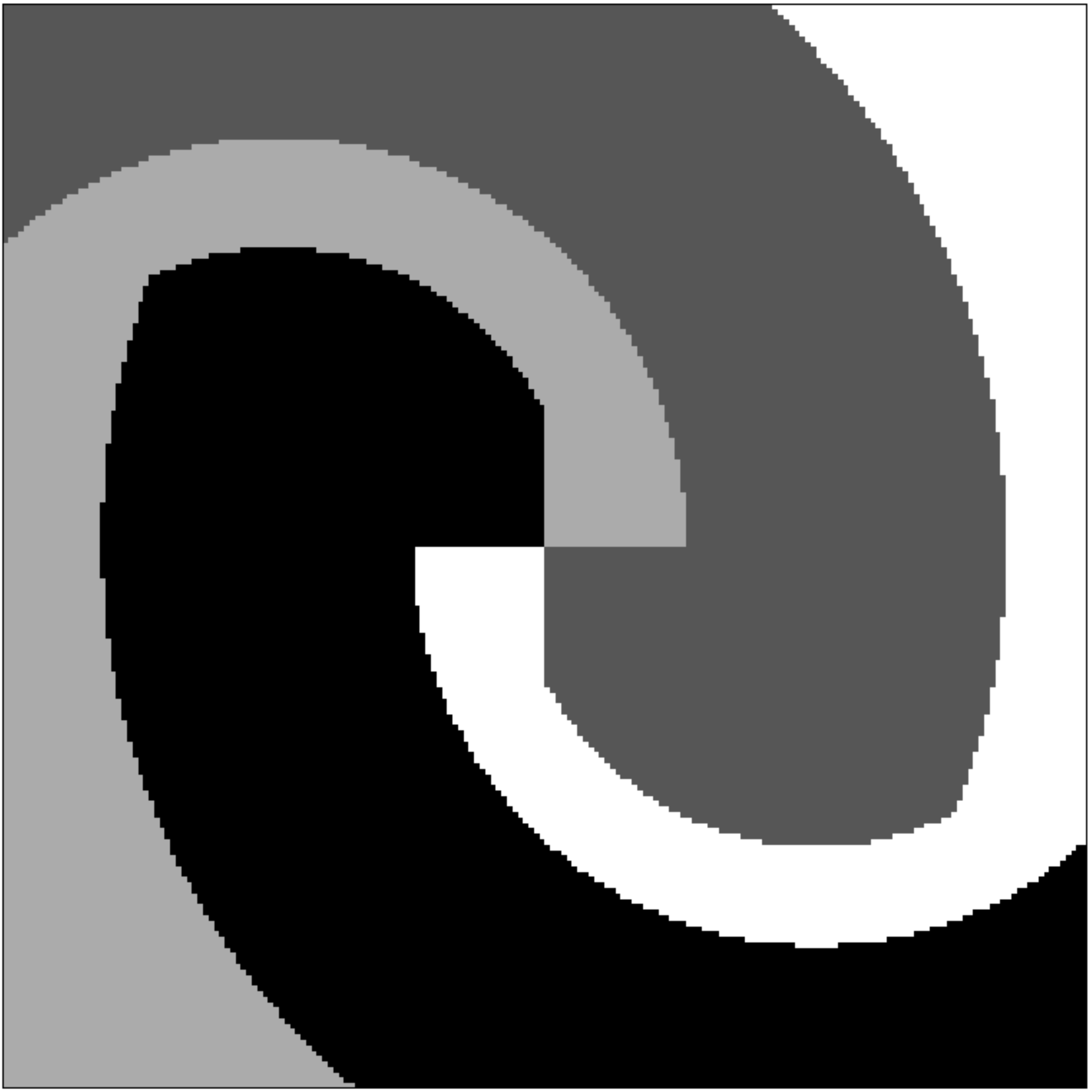} 
  &
  \includegraphics[width=0.4859\textwidth]{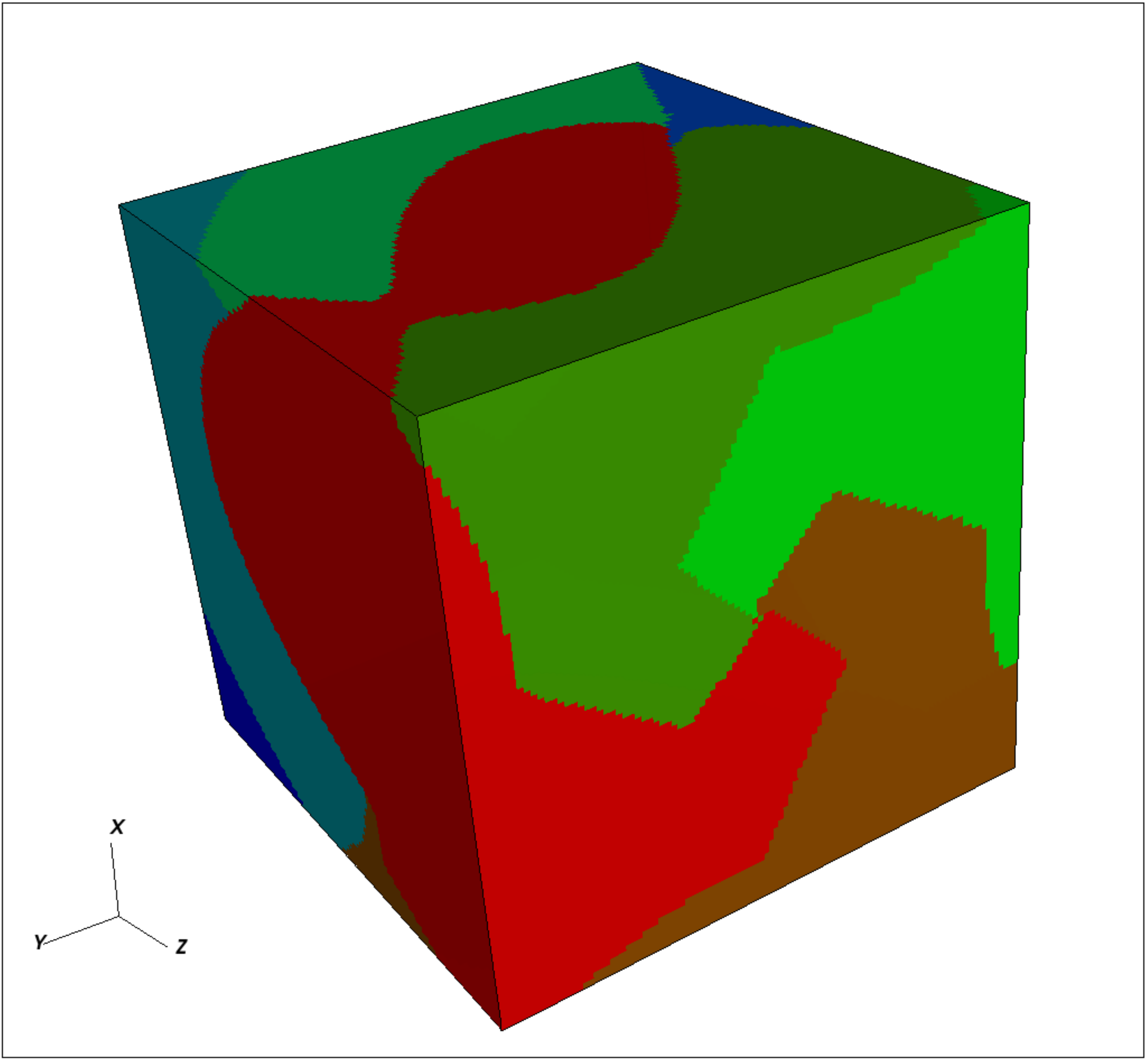}\\
  (a)
  &
  (b)
\end{tabular}
 \caption{Patchy domain decompositions for two test dynamics ($B$ denotes a small ball at the origin): 
 (a) $\Omega=\mathbb{R}^2\setminus B$, $N_p=4$; 
 (b) $\Omega=\mathbb{R}^3\setminus B$, $N_p=8$.}\label{patchy-decomposition}
  \end{center}
\end{figure}

We remark that each patch is a bundle of characteristics enjoying, by construction, the fundamental property of being 
invariant with respect to the optimal dynamics, i.e. $f^*(\Omega_p)\subseteq\Omega_p$. 
This is not completely true, since $f^*$ is built using only the coarse control $\hat a^*$. Moreover, 
at the discrete level, the projection of the patches on the grid $G$ introduces an additional error, in particular if 
the dynamics $f^*$ defines very bended characteristic curves. Figure \ref{patchy-decomposition} shows the patchy decomposition for two test dynamics, 
in dimension two and three respectively. Note that, by construction, the patches do not overlap, sharing only sharp interfaces.

The next step is the parallel computation on the patches. This can be done \emph{avoiding completely} the transmission conditions, 
exploiting the invariance property just discussed above. Since this feature is weakened at the discrete level, it can be 
enforced in the computation by imposing \emph{state constraint} conditions at the interfaces.

Finally, all the solutions are merged together, producing a solution on the whole domain $\Omega$. As expected, this \emph{patchy solution} 
is slightly different from that computed using a standard domain 
decomposition method. Nevertheless, the error is localized at the interfaces between the patches and does not propagate 
in the interior of the sub-domains, 
provided that the grid $G_c$ is not too much coarse compared to the fine grid $G$. This is shown in \cite{CCFP12} by numerical evidence, but 
a precise error estimate is still missing. 

Despite the small errors, the absence of transmission conditions can give to the parallelization 
a remarkable speedup. 
In this respect a key role is played by the relative sizes of the patches. Indeed, we remark that 
the construction of the patchy decomposition is completely driven by the dynamics of the optimal control problem. 
Hence, even a subdivision of the boundary $\partial\Omega$ in sub-sets of the same size can produce 
a highly unbalanced domain decomposition. In these cases the performance of the parallelization is very poor, 
since the processors associated to the smaller patches complete their job (and remain idle) much earlier than 
those corresponding to the larger patches. This drawback was pointed out in \cite{CCFP12} and 
can be overcome via a multi-level technique, which is currently under development. The idea is to alternate 
the construction of the patchy decomposition with the computation of the patchy solution. More precisely, 
one can start and continue the construction of the patches as long as they have about the same size, 
obtaining a first level of balanced sub-domains. The solution is then computed on the first level 
sub-decomposition. The new boundaries (possibly divided again in sub-sets of the same size) are employed 
to start and build the second level sub-decomposition. Moreover, the values of the first level solution 
at the corresponding points (correct values due to the {\em causality property} of hyperbolic equations) 
are assigned as boundary data for the computation of the second level solution. 
This procedure is then iterated and terminates when the sub-decompositions cover the whole domain $\Omega$.

\section{A semi-Lagrangian scheme for second order HJB equations}\label{SL2HJB}
In this section we present a semi-Lagrangian discretization for a special class of stationary second order 
semi-linear equations. 
The resulting scheme will be employed as local solver for the extension of the patchy method to this more general setting, 
as discussed in the next section. 

The prototype problem we have in mind is the following boundary value problem for second order semi-linear equations, written in a control theory perspective:
\begin{equation}\label{HJB2}
\left\{ \begin{array}{ll}
         \displaystyle\max_{a\in A}\left\{ -\mathcal{D}(x,a)u(x) -l(x,a) \right\} = 0 & x\in \Omega\\
         u(x)=g(x) & x\in \partial \Omega
        \end{array}
\right.
\end{equation}
where $\Omega\subset\mathbb{R}^n$ is an open set, $u:\Omega\to\mathbb{R}$, $A\subset\mathbb{R}^m$ is a compact set and 
the second order differential operator 
$\mathcal{D}$ is given by
$$
\mathcal{D}(x,a)=\frac12\sum_{i,j=1}^n\left(\sum_{k=1}^d \sigma_{ik}(x,a)\sigma_{jk}(x,a)\right)\frac{\partial^2}{\partial x_i \partial x_j} + 
\sum_{i=1}^n f_i (x,a) \frac{\partial}{\partial x_i}
$$
with $\sigma : \Omega\times A\to\mathcal{L}(\mathbb{R}^d;\mathbb{R}^n)$, $f:\Omega\times A\to \mathbb{R}^n$ and 
$l:\Omega\times A\to \mathbb{R}$, $g:\partial\Omega\to \mathbb{R}$.

The crucial point here is that the connection with deterministic optimal control, depicted in the previous section, 
is lost due to the presence of the diffusion term $\sigma$. Indeed, in this case, characteristic curves are no longer well defined by 
the system of controlled ordinary differential equations $(\ref{ODE})$. 
Nevertheless, a weak notion of characteristics is still available, 
interpreting these curves as \emph{generalized trajectories}, namely solutions to the following system of controlled \emph{stochastic} 
differential equations with dynamics $f$:
\begin{equation}\label{SDE}
\left\{ \begin{array}{ll}
         dX(t) = f(X(t),\alpha(t))dt +\sigma(X(t),\alpha(t)) dW(t)\,, & t>0\\
         X(0)  = x & 
        \end{array}
\right. 
\end{equation}
Here $X(t)$ is a progressively measurable process, representing the state of a system evolving in $\Omega$ starting from 
$x$, the process 
$\alpha(t)$ is the control applied to the system at the time $t$ with values in the control set $A$ and $W(t)$ is a $d$-dimensional 
Wiener process. We define the set of admissible controls 
$$\mathcal{A}=\{\alpha :[0,+\infty)\to A,\mbox{ progressively measurable}\}$$
and we denote by $X(t;x,\alpha(\cdot))$ the 
solution to $(\ref{SDE})$ starting from $x$ using the control $\alpha(\cdot)$. Finally, we define the first time arrival to the boundary $\partial \Omega$ as
$$
\tau(x,\alpha(\cdot))=\inf\left\{t\geq 0\,:\, X(t;x,\alpha(\cdot))\in\partial\Omega\right\} 
$$
and we consider the following cost functional ($\mathbb{E}$ stands for the probabilistic expectation):
$$
J\big(x,\alpha(\cdot)\big)=\mathbb{E}\left\{\int_0^{\tau(x,\alpha(\cdot))} l\big(X(s;x,\alpha(\cdot)),\alpha(s)\big)\,ds + g\big(X(\tau(x,\alpha(\cdot)))\big)\right\}\,.
$$
In this setting, the unique viscosity solution $u$ to the problem $(\ref{HJB2})$ 
can be interpreted as the value function of the following \emph{stochastic} minimum time problem with 
running cost $l$ and exit cost $g$:
$$
u(x)=\inf_{\alpha(\cdot)\in\mathcal{A}} J\big(x,\alpha(\cdot)\big) \qquad x\in\Omega\,.
$$
Typical assumptions on $f$, $g$, $l$ and $\sigma$ for the well posedness of the problem is boundedness and Lipschitz continuity in space uniformly 
in the control. We refer the interested reader to  \cite{CF95} and the references therein for further explanations and details.

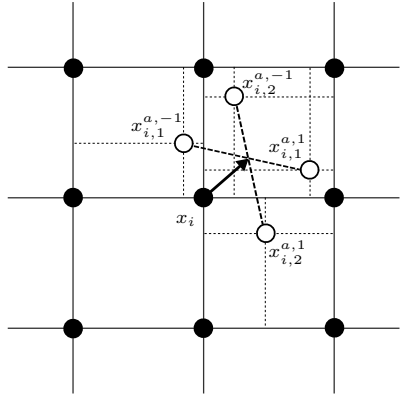
\begin{figure}[htp]
\begin{center}
 \includegraphics[width=.4\textwidth]{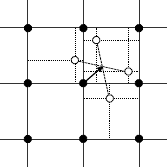}\\
  \caption{An example in dimension $2$ with $\sigma:\Omega\times A\to\mathcal{L}(\mathbb{R}^2;\mathbb{R}^2)$. 
  The arrow represents the vector field $f(x_i,a)$ for a generic control $a\in A$.  For $s=\pm 1$ and $k=1,2$ the values at (white) points $x_{i,k}^{a,s}$
  are reconstructed via bilinear interpolation of (black) grid nodes.}\label{bilinear}
  \begin{picture}(1,1)
  \put(-10,115){$\scriptstyle x_i$}
  \put(-27,150){$\scriptstyle x_{i,1}^{a,-1}$}
  \put(25,142){$\scriptstyle x_{i,1}^{a,1}$}
   \put(15,166){$\scriptstyle x_{i,2}^{a,-1}$}
   \put(25,102){$\scriptstyle x_{i,2}^{a,1}$}
   \end{picture}
  \end{center}
\end{figure}

Following \cite{CF95}, the semi-Lagrangian discretization of equation $(\ref{HJB2})$ can be performed introducing a time step $h$, 
interpreting the first order term in the operator $\mathcal{D}$ as a directional derivative of $u$ along the dynamics $f$ and 
approximating the Wiener process in the stochastic differential equation $(\ref{SDE})$ by means of discrete time 
Markov chains. 
We introduce in the state space a structured grid $G$ with uniform step $\Delta x$ in each coordinate direction and nodes $x_i$ 
for $i=1,...,N$. This leads to the following scheme, which is a nonlinear system in the form of a fixed point operator:
\begin{equation}\label{SLscheme2}
\left\{
\begin{array}{ll}
x_{i,k}^{a,s}=x_i+hf(x_i,a)+s\sqrt{h}\sigma_k (x_i,a) & s=\pm 1,\, k=1,...,d\\
u(x_i) = \displaystyle\min\limits_{a\in A}\left\{ \frac{1}{2d} \sum_{k=1}^d \sum_{s=\pm 1} u(x_{i,k}^{a,s})+ h\,l(x_i,a)\right\}
& x_i \in G\cap \Omega\\
u(x_i)=g(x_i) &  x_i \in G\cap \partial \Omega
\end{array}
\right.
\end{equation}
where $\sigma_k$ denotes the $k$-th row of $\sigma$. The evaluation of the $\min$ operator is done by direct comparison 
discretizing the control set $A$ with $N_A$ points. Moreover, as usual in semi-Lagrangian approximations, the $2d$ points $x_{i,k}^{a,s}$ 
do not lie in general on the grid and the values of $u$ at the these points need to be reconstructed by interpolation. 
A natural choice here is the bilinear interpolation, as shown in Figure \ref{bilinear}.

We choose a time step $h$ such that, for each node $x_i\in G$ and every control $a\in A$, the points $x_{i,k}^{a,s}$ fall in the first 
neighboring cells of $x_i$. The reason is twofold. First, it is easy to prove that 
the semi-Lagrangian scheme $(\ref{SLscheme2})$ is consistent to equation $(\ref{HJB2})$ only with order 
$\mathcal{O}(h)$ 
and we do not want to get a too poor accuracy in the approximation of the trajectories. Second and more important, 
we want to keep the stencil of the scheme strictly localized, in order to reduce the computational efforts.

Now, let us denote by 
$x_{i_0,k}^{a,s}$, $x_{i_1,k}^{a,s}$, $x_{i_2,k}^{a,s}$, $x_{i_3,k}^{a,s}$ the four points involved in the reconstruction of 
the value 
$u(x_{i,k}^{a,s})$, with weights $\lambda_{i_0,k}^{a,s}$, $\lambda_{i_1,k}^{a,s}$, $\lambda_{i_2,k}^{a,s}$, 
$\lambda_{i_3,k}^{a,s}$ respectively. Without loss of generality, we can assume that $x_{i_0,k}^{a,s}=x_i$ for each 
$k=1,...,d$ and $s=\pm 1$. It follows that
$$
u(x_{i,k}^{a,s})=\sum_{p=0}^3 \lambda_{i_p,k}^{a,s} u(x_{i_p,k}^{a,s})=
\lambda_{i_0,k}^{a,s} u(x_i)+\sum_{p=1}^3 \lambda_{i_p,k}^{a,s} u(x_{i_p,k}^{a,s})\,.
$$
Substituting this expression in the scheme $(\ref{SLscheme2})$ we get 
\begin{equation}\label{self-dependency}
u(x_i) = \displaystyle\min\limits_{a\in A}\left\{ \frac{1}{2d} \sum_{k=1}^d \sum_{s=\pm 1} \lambda_{i_0,k}^{a,s} u(x_i)+
\frac{1}{2d} \sum_{k=1}^d \sum_{s=\pm 1}\sum_{p=1}^3 \lambda_{i_p,k}^{a,s} u(x_{i_p,k}^{a,s}) + h\,l(x_i,a)\right\}\,,
\end{equation}
which shows explicitly the dependency of $u(x_i)$ on itself. 

As explained in the 
Introduction, in the hyperbolic case ($\sigma\equiv 0$) relevant information starts from the boundary of the domain, 
and propagates along characteristics backward in time at finite speed. At the discrete level, we can try to mimic 
this \emph{causality property}, by means of an appropriate ordering of the grid nodes that tracks the front of information. 
In this way, we can get a cascade effect that decouples, at least partially, the nonlinear system, 
drastically reducing the number of iterations to reach convergence. 
Despite this property does not hold in the continuous case in presence of diffusion ($\sigma\not\equiv 0$), we can still have a 
\emph{numerical causality} if some suitable condition on the parameters is satisfied. We will come back on this important point 
in the next section.

Here the key argument is that we have to remove in $(\ref{self-dependency})$ the self-dependency on $u(x_i)$, 
which makes the scheme strongly iterative by construction. 
This is a really crucial step, if we hope to benefit of the speedup induced by the causality property.  
To this end, we denote by 
$$
\Lambda_0(a)=\frac{1}{2d} \sum_{k=1}^d \sum_{s=\pm 1} \lambda_{i_0,k}^{a,s}\,,\qquad 
\Lambda_1(a)=\frac{1}{2d} \sum_{k=1}^d \sum_{s=\pm 1}\sum_{p=1}^3 \lambda_{i_p,k}^{a,s} u(x_{i_p,k}^{a,s}) + h\,l(x_i,a)\,,
$$
so that
\begin{equation}\label{SLimplicit}
 u(x_i) = \displaystyle\min\limits_{a\in A}\left\{  \Lambda_0(a)u(x_i)+\Lambda_1(a)\right\}\,.
\end{equation}
We remark that $\Lambda_0(a)$ and $\Lambda_1(a)$ do not depend on the value $u(x_i)$. Moreover, we note that 
$\Lambda_0(a)<1$ for all $a\in A$, since the interpolation weights $\lambda_{i_0,k}^{a,s}$ can be simultaneously equal to $1$ only 
if the time step $h=0$. Then, it is well defined the value
\begin{equation}\label{SLexplicit}
 v(x_i) = \displaystyle\min\limits_{\bar a\in A}\left\{ \frac{\Lambda_1(\bar a)}{1-\Lambda_0(\bar a)}\right\}
\end{equation}
and we claim that $v(x_i)=u(x_i)$. 
This kind of explicitation is trivial in the linear case without diffusion (i.e. for $f$ not depending 
on the control $a$ and $\sigma\equiv 0$), but it is less evident (and to our knowledge also quite surprising) 
in the general nonlinear case. 
So let $a^*\in A$ be the control achieving the minimum in $(\ref{SLimplicit})$, i.e. 
$$
u(x_i) = \Lambda_0(a^*)u(x_i)+\Lambda_1(a^*)\,.
$$
This implies
$$
u(x_i)=\frac{\Lambda_1(a^*)}{1-\Lambda_0(a^*)}\geq 
\displaystyle\min\limits_{\bar a\in A}\left\{ \frac{\Lambda_1(\bar a)}{1-\Lambda_0(\bar a)}\right\}=v(x_i)\,.
$$
Conversely, let $\bar a^*\in A$ be the control achieving the minimum in $(\ref{SLexplicit})$, i.e. 
$$
v(x_i) = \frac{\Lambda_1(\bar a^*)}{1-\Lambda_0(\bar a^*)}\,.
$$
This implies  
$$
u(x_i) \leq  \Lambda_0(\bar a^*)u(x_i)+\Lambda_1(\bar a^*)\,.
$$
Since $\Lambda_0(a^*)<1$, we get 
$$
u(x_i)\leq \frac{\Lambda_1(a^*)}{1-\Lambda_0(a^*)}=v(x_i)\,.
$$
Finally, we obtain the following modified scheme
\begin{equation}\label{SLscheme3}
\left\{
\begin{array}{ll}
x_{i,k}^{a,s}=x_i+hf(x_i,a)+s\sqrt{h}\sigma_k (x_i,a) & \hskip-1cm s=\pm 1,\, k=1,...,d\\ & \\
u(x_i) = \displaystyle\min\limits_{a\in A}
       \left\{
\frac{ 
      \displaystyle\frac{1}{2d} \sum_{k=1}^d \sum_{s=\pm 1}\sum_{p=1}^3 \lambda_{i_p,k}^{a,s} u(x_{i_p,k}^{a,s}) + h\,l(x_i,a)
     }
     {
      1-\displaystyle\frac{1}{2d} \sum_{k=1}^d \sum_{s=\pm 1} \lambda_{i_0,k}^{a,s}
     }
       \right\}
& x_i \in G\cap \Omega\\ & \\
u(x_i)=g(x_i) &  x_i \in G\cap \partial \Omega
\end{array}
\right.
\end{equation}
where the value at each node $x_i$ now depends only on values at nodes different from $x_i$. 
In Section \ref{experiments} we compare this 
scheme with the original scheme $(\ref{SLscheme2})$, showing the effectiveness of the modification in terms of iterations to reach 
convergence.

We conclude this section by remarking that the proposed semi-Lagrangian scheme can handle by construction 
very degenerate diffusions. Indeed, at points where the diffusion coefficient $\sigma=0$, 
we naturally recover the semi-Lagrangian scheme presented in \cite{CCFP12} 
for first order Hamilton-Jacobi-Bellman equations (see also the numerical experiments in Section \ref{experiments}). 
This is not the case for other types of discretization, e.g. finite elements, 
where the degeneracy in the diffusion produces instabilities that need to be treated with very specific techniques. 
This is due to the fact that the variational forms associated to the equations under consideration suffer a lack of coercivity. 
We refer the interested reader to \cite{QV99} for details and insights on this topic.  

\section{Patchy decomposition for second order semi-linear equations}\label{PD2}
In this section we aim to extend the patchy decomposition method to the class of second order semi-linear equations 
presented in the previous section.

The main idea is really simple: we first compute the patchy domain decomposition in the hyperbolic case 
with $\sigma\equiv 0$. Then we solve in parallel the full equation with $\sigma\not\equiv 0$ using the patchy decomposition.

At a first look this attempt could seem meaningless. Indeed, due to the diffusion, information 
spreads instantaneously from the boundary to the whole domain and then it cannot be confined in independent sub-domains, as for the 
first order case. This implies that, in order to compute the correct solution in this more general setting, 
we cannot replace the transmissions between the patches with state constraint boundary conditions, as discussed in Section \ref{PD1}. 
So we loose the main advantage of the patchy decomposition compared to an arbitrary and static decomposition. 
Moreover, we recall 
that the preliminary step in the computation of the dynamic decomposition results in an additional cost in terms of CPU time.
The reader can easily convince himself that, in presence of transmission conditions, the performance of a parallel computation 
does not significantly depend, in general, on the particular shape of the sub-domains, but only on their number. 
As pointed out in \cite{CCFP12}, if we assume that 
the sub-domains have about the same size, then the overall time consumption to compute the solution almost exclusively depends 
on the number of processors involved in the computation and the transmission delays.

Nevertheless, we show through the next sections that the patchy method with transmission conditions can still be competitive 
if we combine two different features, described in the following sub-sections.
\subsection{The upwind diffusion ball condition}
Here we present the first ingredient for the extension of the patchy method to second order semi-linear equations. 
In particular we show that, under suitable relations between the parameters, the discretization of equation $(\ref{HJB2})$ 
behaves in a sense more like hyperbolic 
than elliptic.

To simplify the presentation we consider a special case in dimension two. Let $\varepsilon>0$ and take 
$\sigma: \Omega\to\mathcal{L}(\mathbb{R}^2;\mathbb{R}^2)$ of the form
\begin{equation}\label{const-eps}
\sigma=\sqrt{2\varepsilon}\left(\begin{array}{cc}
                                 1 & 0 \\ 0 & 1
                                \end{array}
\right)
\end{equation}
Equation $(\ref{HJB2})$ writes
$$
\left\{ \begin{array}{ll}
         -\varepsilon\Delta u(x)+\displaystyle\max_{a\in A}\left\{ -f(x,a)\cdot \nabla u(x) -l(x,a) \right\} = 0 & x\in \Omega\\
         u(x)=g(x) & x\in \partial \Omega
        \end{array}
\right.
$$
and the semi-Lagrangian scheme $(\ref{SLscheme3})$ simplifies in 
$$
\left\{
\begin{array}{ll}
x_{i,k}^{a,s}=x_i+hf(x_i,a)+s\sqrt{2\varepsilon h}\, e_k & \hskip-1cm s=\pm 1,\, k=1,2\\ & \\
u(x_i) = \displaystyle\min\limits_{a\in A}
       \left\{
\frac{ 
      \displaystyle\frac{1}{4} \sum_{k=1}^2 \sum_{s=\pm 1}\sum_{p=1}^3 \lambda_{i_p,k}^{a,s} u(x_{i_p,k}^{a,s}) + h\,l(x_i,a)
     }
     {
      1-\displaystyle\frac{1}{4} \sum_{k=1}^2 \sum_{s=\pm 1} \lambda_{i_0,k}^{a,s}
     }
       \right\}
& x_i \in G\cap \Omega\\ & \\
u(x_i)=g(x_i) &  x_i \in G\cap \partial \Omega
\end{array}
\right.
$$
where $e_k$ denotes the $k$-th element of the canonical base of $\mathbb{R}^2$. 

As discussed in Section \ref{SL2HJB}, 
for each node $x_i$ and control $a\in A$, 
we want to reconstruct the value of $u$ at the four points $x_{i,k}^{a,s}$ ($k=1,2$, $s=\pm 1$), 
via bilinear interpolation of the first neighboring grid nodes. To this end, it suffices to choose a time step $h$ such that 
\begin{equation}\label{upper-bound}
hf_{\max}+\sqrt{2\varepsilon h}<\Delta x\,,
\end{equation}
where $f_{\max}$ denotes the maximum of the dynamics $f$ on $\overline{\Omega}\times A$. 
Note that condition $(\ref{upper-bound})$ can be localized for each node $x_i$ and control $a\in A$. This results in 
a more accurate approximation, but it is clearly more expensive from the computational point of view. 

\begin{figure}[htp]
\begin{center}
 \includegraphics[width=0.85\textwidth]{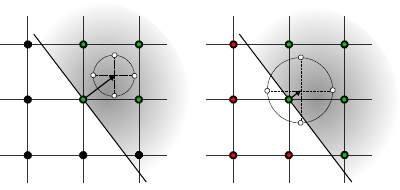}\\
 \begin{tabular}{ccc}
  \hspace{-.75cm}(a) &\hspace{4.5cm} & (b)
 \end{tabular}
 \caption{The upwind diffusion ball condition: (a) $B\subseteq \pi$, the interpolating (green) grid nodes are upwind w.r.t. the vector field $f$; 
 (b) $B\not\subseteq \pi$, some downwind (red) grid nodes enter in the computation.}\label{halfspace}
  \end{center}
\end{figure}

Now consider the half space defined by $\pi=\{ x\in\mathbb{R}^2\,:\, f(x_i,a)\cdot (x-x_i)\geq 0\}$ and also the ball 
$B$ of radius $\sqrt{2\varepsilon h}$ centered at 
$x_i+hf(x_i,a)$, enclosing the four points that need to be reconstructed (see Figure \ref{halfspace}). 
The key argument of our construction 
is that, if the \emph{diffusion ball} $B$ is contained in the half space $\pi$, then the value $u(x_i)$ is computed using only 
grid nodes that are \emph{upwind} with respect to the vector field $f$, as shown in Figure \ref{halfspace}-(a). 
If we combine this property with a suitable order for processing the grid nodes (aka causality, see next subsection) we can accelerate the 
convergence of the scheme in the same spirit of fast-marching methods, reducing significantly the number of iterations. 
On the contrary, if part of $B$ crosses $\pi$, then also downwind nodes are employed in the reconstruction of 
$u(x_i)$, as shown in Figure \ref{halfspace}-(b). It follows that some information is flowing in directions 
opposite to the vector field $f$, so that $u(x_i)$ cannot be computed in a single pass fashion and 
the number of iterations to reach convergence increases. This is clearly 
expected, since we are considering the particular example $(\ref{const-eps})$ in which the diffusion uniformly spreads information 
in all the directions. 
But the crucial point here is that the semi-Lagrangian scheme propagates diffusion at speed $\sqrt{2\varepsilon}$, 
and not instantaneously 
as in the continuous case. Indeed, at each iteration, we move with discrete time steps from the point $x_i$ to 
the points $x_{i,k}^{a,s}$ ($k=1,2$, $s=\pm 1$), adding the contribution of the vector field $f$ and the diffusion.  
Then, the upwind condition on the diffusion ball holds if the time step $h$ satisfies also 
\begin{equation}\label{lower-bound}
hf_{\min}-\sqrt{2\varepsilon h}>0\,,
\end{equation}
where $f_{\min}$ denotes the minimum of the dynamics $f$ on $\overline{\Omega}\times A$.
Note that also \eqref{lower-bound} can be localized for each node $x_i$ and control $a\in A$, but in the following we will 
always consider only global conditions.

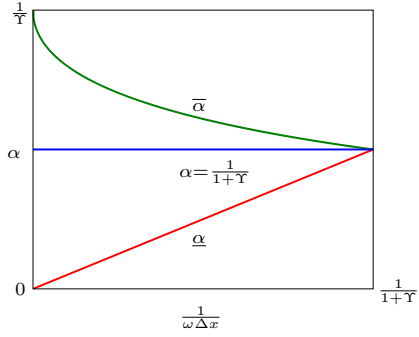
\begin{figure}[htp]
\begin{center}
 \includegraphics[width=0.4\textwidth]{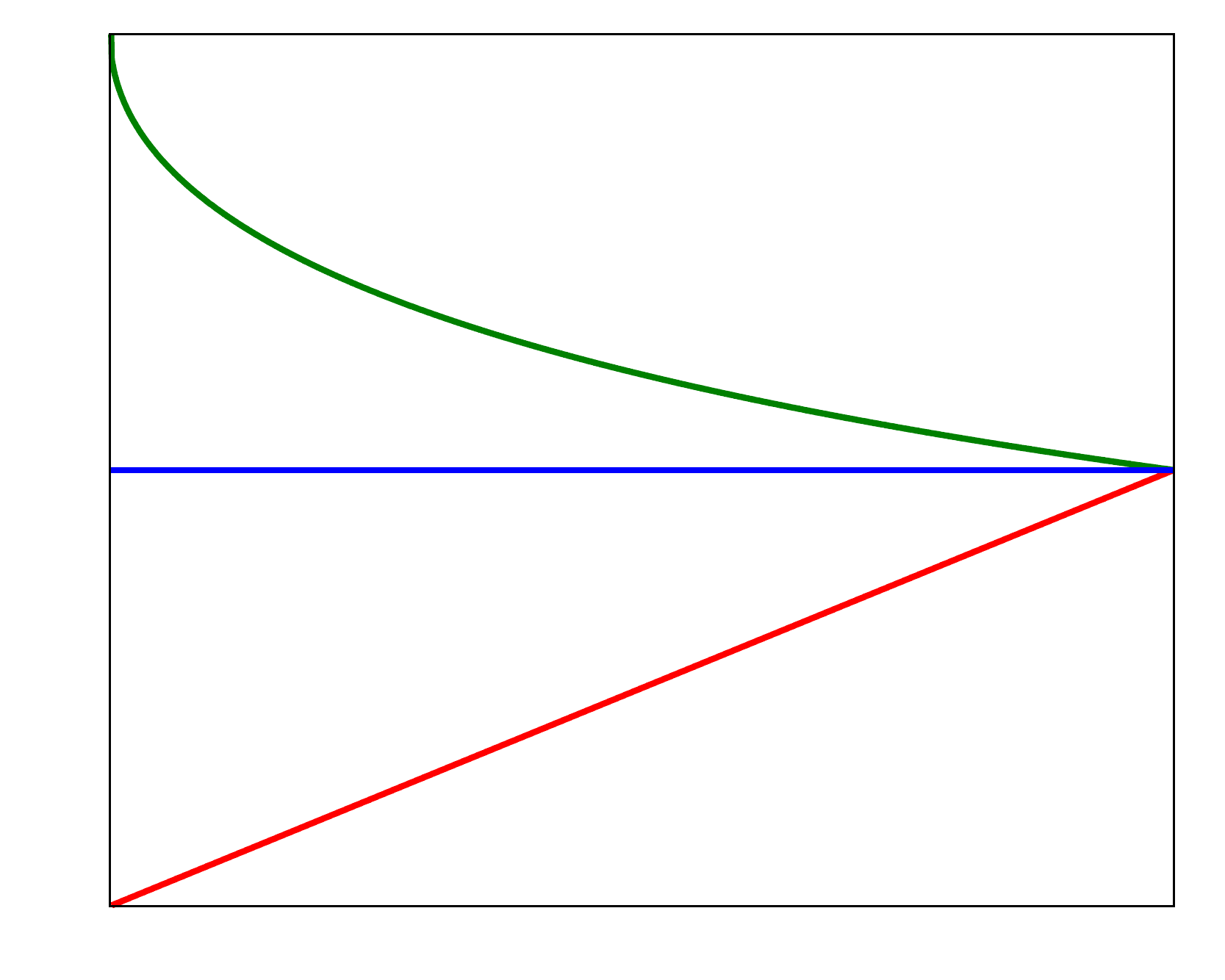}\\
 \caption{The value $\alpha=\frac{1}{1+\Upsilon}$ (in blue) and the graphs of $\underline{\alpha}$, 
 $\overline{\alpha}$ (in red and green respectively) as functions of $1/\omega\Delta x$,
 ranging in $(0,\frac{1}{1+\Upsilon})$ as prescribed by condition $(\ref{hyperbolicity})$.}\label{alpha-opt}
   \begin{picture}(1,1)
  \put(-5,95){$\scriptstyle\alpha=\frac{1}{1+\Upsilon}$}
  \put(0,70){$\scriptstyle \underline{\alpha}$}
  \put(0,120){$\scriptstyle \overline{\alpha}$}
  \put(-70,102){$\scriptstyle \alpha$}
  \put(-5,40){$\scriptstyle \frac{1}{\omega\Delta x}$}
  \put(-67,50){$\scriptstyle 0$}
   \put(-69,155){$\scriptstyle \frac{1}{\Upsilon}$}
  \put(70,50){$\scriptstyle \frac{1}{1+\Upsilon}$}
   \end{picture}
  \end{center}
\end{figure}

We now relate the discretization parameters to the data, aiming to satisfy both conditions $(\ref{upper-bound})$ and 
$(\ref{lower-bound})$. To this end, we set $h=\alpha\Delta x/f_{\min}$ with $\alpha>0$ to be determined. Moreover,  
we define $\omega=f_{\min}/2\varepsilon$ and $\Upsilon=f_{\max}/f_{\min}$. Substituting in $(\ref{upper-bound})$, we get
$$
\sqrt{\frac{\alpha\Delta x}{\omega}}<(1-\alpha\Upsilon)\Delta x\,,
$$
which implies, for $\alpha<1/\Upsilon$ and by straightforward computations, the following condition:
\begin{equation}\label{upper-alpha}
 \alpha< \frac{1}{\Upsilon}-\frac{\sqrt{1+4\omega\Delta x\Upsilon}-1}{2\omega\Delta x\Upsilon^2}=:\overline{\alpha}
\end{equation}
with $\overline{\alpha}$ satisfying $0<\overline{\alpha}<1/\Upsilon$. On the other hand, substituting in $(\ref{lower-bound})$, we immediately get
\begin{equation}\label{lower-alpha}
 \alpha> \frac{1}{\omega\Delta x}=:\underline{\alpha}\,.
\end{equation}
By imposing the compatibility of the two conditions, i.e. $\underline{\alpha}<\overline{\alpha}$, we easily obtain the following 
\emph{upwind diffusion ball condition}:
\begin{equation}\label{hyperbolicity}
\frac{1}{\omega}<\frac{\Delta x}{1+\Upsilon}\,.
\end{equation}
Summarizing, if condition $(\ref{hyperbolicity})$ is satisfied, we can choose $\alpha\in(\underline{\alpha},
\overline{\alpha})$ such that both $(\ref{upper-bound})$ and 
$(\ref{lower-bound})$ hold. In particular, it is easy to check that $\alpha=\frac{1}{1+\Upsilon}$ 
is the only value satisfying $\underline{\alpha}<\alpha<\overline{\alpha}$ for all $\omega$ and $\Delta x$ such that $(\ref{hyperbolicity})$ holds 
(see also Figure \ref{alpha-opt}). 
Conversely, if $(\ref{hyperbolicity})$ fails, only one of the two conditions can be fulfilled, by choosing $\alpha$ appropriately. 
We stress again that, for accuracy and computational issues discussed in Section \ref{SL2HJB}, 
from now on we will always satisfy condition $(\ref{upper-bound})$, possibly loosing condition $(\ref{lower-bound})$. 

In the general case $(\ref{HJB2})$, where the diffusion coefficient is given by the matrix-valued map 
$\sigma : \Omega\times A\to\mathcal{L}(\mathbb{R}^d;\mathbb{R}^n)$, we obtain the same condition $(\ref{hyperbolicity})$ with 
$\omega=f_{\min}/\|\sigma\|_\infty^2$, where 
$$
\|\sigma\|_\infty=\adjustlimits\max_{(x,a)\in\overline{\Omega}\times A}\max_{k=1,..,d}\sqrt{\sum_{m=1}^n \sigma_{km}^2(x,a)}\,.
$$
We remark that $\|\sigma\|_\infty$ identifies the diffusion vector with greatest length among the rows of $\sigma$, namely an upper bound for the radius of the 
diffusion ball.

For a fixed and small mesh size $\Delta x$, condition $(\ref{hyperbolicity})$ can be clearly fulfilled if 
the controlled advection \emph{dominates} diffusion, i.e. for $\omega>>1$. 
Note that $\omega$ is a characteristic quantity of the problem. It resembles the Reynolds number for Navier-Stokes 
equations and the case $\omega>>1$ is of great interest in the applications \cite{QV99}.
In this perspective, condition $(\ref{hyperbolicity})$ 
describes an interesting \emph{threshold effect}, by means of the mesh dependent 
parameter $\tau_{\Delta x}=\Delta x/(1+\Upsilon)$. Indeed, for $1/\omega<\tau_{\Delta x}$, the upwind diffusion ball condition is satisfied, and 
the semi-Lagrangian scheme exhibits the same peculiarities of the 
hyperbolic case. Conversely, for $1/\omega>\tau_{\Delta x}$, 
the upwind diffusion ball condition fails and an elliptic behavior emerges. 
Note that this effect is also related to the degree of anisotropy $\Upsilon\ge 1$. 
Equivalently, we can fix $\omega$, depending on the data $f$, $\sigma$, 
and choose the mesh size $\Delta x$ according 
to the threshold $\tau_\omega=(1+\Upsilon)/\omega$. 
It follows that, at a coarse scale, i.e. for $\Delta x>\tau_\omega$, the approximation 
``looks like'' that of a controlled advection equation. Again, this effect also depends on the degree of anisotropy, 
in the sense that the larger is $\Upsilon$ the coarser should be the mesh. 
On the other hand, at a sufficiently fine scale $\Delta x<\tau_\omega$, 
diffusion starts to get noticed and this hyperbolic behavior disappears as $\Delta x\to 0$. 
This is expected, due to the consistency of the semi-Lagrangian scheme with the considered semi-linear equation.

\begin{figure}[htp]
\begin{center}
 \includegraphics[width=0.4\textwidth]{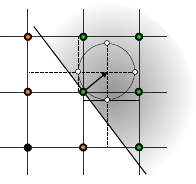}\\
 \caption{The upwind diffusion ball condition is satisfied, but some downwind (orange) grid nodes enter in the computation with small 
 interpolation weights.}\label{halfspace-ex}
  \end{center}
\end{figure}

Finally, we remark that, from a numerical point of view, the effect of the upwind diffusion ball condition 
$(\ref{hyperbolicity})$ on the performance of the scheme can be less relevant than expected. Indeed, 
even if the diffusion ball $B$ is contained in the half space $\pi$ associated to $f$, 
some grid points involved in the 
interpolation can be outside, as shown in Figure \ref{halfspace-ex}. 
This strongly depends on the alignment to the grid of the optimal vector field 
achieving the minimum in $(\ref{SLscheme3})$. 
Nevertheless, as confirmed by the numerical experiments in Section \ref{experiments}, this drawback does not significantly affect the hyperbolic behavior of 
the scheme, since the downwind grid points contribute in these cases with very small interpolation weights, 
compared to the upwind ones. 
\subsection{The role of causality}
Under the regime given by condition $(\ref{hyperbolicity})$, we have a chance to get a good performance of the 
patchy method. To this end, the second crucial point for letting work all the machinery in this more general setting 
is the causality property of hyperbolic equations outlined in the previous section, 
a quite well established topic in the literature. 

Starting from the work by Tsitsiklis \cite{T95} and then Sethian \cite{S96}, a lot of research has 
been devoted to find an implicit order of the nodes of a grid that allows to compute the solution in just one or very few 
iterations, the so called \emph{single-pass property}. With this idea in mind, the celebrated fast-marching method has been proposed to solve the eikonal 
equation. It can be proved that, in this special case, 
the right order corresponds to process the nodes by ascending values, progressively accepting as correct 
(and then removing from the computation) the node with the minimal value. This translates in computing the 
solution following its level sets, namely propagating information along its gradient curves. 
Since for the eikonal equation gradient curves coincide with characteristics, we get the correct solution. 

Unfortunately, this ordering is optimal only for eikonal equations of isotropic type, 
namely for a dynamics $f$ that does not depend (or depends in a weak sense) on the control variable. 
We refer the reader to \cite{SV03} for a detailed 
explanation of this problem. Whenever a strong anisotropy comes into play, the choice of a correct order for the grid nodes 
becomes a subtle topic, 
meaning that it may not even exist, despite the causality property always holds at the continuous level. 
This goes beyond the scope of this paper, but we refer the reader to \cite{CCF14a} for a discussion on the applicability of 
fast-marching-like methods to general Hamilton-Jacobi equations. We only point out that, in order to compute the correct solution, 
one has to enlarge the stencil of the scheme and/or give up the single-pass property mentioned above. This results in more 
computational efforts and more iterations for the scheme to reach convergence.

We mention here also the fast-sweeping method \cite{TCOZ04,Z05} and some of its generalizations  \cite{CCF14}, based on another technique 
that, in a weaker sense, still exploits the causality property. The grid nodes are alternatively swept in a pre-determined 
number of directions according to the dimension of the problem, until convergence is reached. 
This makes the method iterative by construction, but it allows to compute the solution to very general equations of Hamilton-Jacobi type. 
The number of iterations to reach convergence is strongly dependent on the problem and the mesh structure. 
Nevertheless, it can be proved that only $2n$ sweeps are needed to compute the solution to the eikonal equation in dimension $n$.

In Section \ref{PD1}, we denoted by $\hat u$ the interpolation on the fine grid $G$ of the coarse solution $u_c$ on the coarse grid $G_c$. 
Our strategy is then to sort the grid nodes of $G$ in a fast-marching fashion, according to the increasing values of $\hat u$. 
This can be accomplished with some additional but cheap time consumption, using some state-of-the-art sorting algorithm, embedded in the pre-computation step.
Note that this procedure was already mentioned in \cite{CCFP12} as a possible add-on for the original patchy method. 
Here it becomes an essential part of the proposed method. We will see in the next section that this ordering of the grid nodes, 
even if in general sub-optimal, can give an exceptional speedup to the computation. 

\subsection{The patchy algorithm}
Here we summarize all the implementation steps of the patchy method for second order semi-linear equations.  
To this end, we present both the discrete version of the procedure $(\ref{synthesis})$ 
for the synthesis of the feedback optimal control and the atcual construction of the patchy decomposition discussed in Section \ref{PD1}.

We consider the \emph{non modified} semi-Lagrangian scheme $(\ref{SLscheme2})$ in the special case without diffusion, 
i.e. $\sigma\equiv 0$:
$$
\left\{
\begin{array}{ll}
u(x_i) = \displaystyle\min\limits_{a\in A}\left\{ u(x_i+hf(x_i,a))+ h\,l(x_i,a)\right\}
& x_i \in G\cap \Omega\\
u(x_i)=g(x_i) &  x_i \in G\cap \partial \Omega
\end{array}
\right.
$$
Once the value function $u$ is computed, we easily get
\begin{equation}\label{discrete-synthesis}
\begin{array}{ll}
a^*(x_i) = \displaystyle\arg\min\limits_{a\in A}\left\{ u(x_i+hf(x_i,a))+ h\,l(x_i,a)\right\}
& x_i \in G\cap \Omega\,.
\end{array}
\end{equation}
Note that the computational cost of this operation mainly depends on the number of points $N_A$ used to discretize 
the control set $A$, but also on the bilinear interpolation for the reconstruction of $u$ at the points $x_i+hf(x_i,a)$.

We proceed with the construction of the dynamic decomposition. We divide the boundary nodes of $G$ in a fixed number $N_P$ of 
sub-sets $\Gamma_p$, with $p=1,...,N_P$, defining the initial guess 
$$
\phi_p(x_i)=\left\{ \begin{array}{ll}
         1  & x_i\in\Gamma_p\\
         0 & \mbox{otherwise}  \end{array}
\right. 
$$
Then, for each $p$ and a fixed tolerance $\varepsilon_P$, we iterate until convergence the scheme
\begin{equation}\label{discrete-advection}
 \begin{array}{ll}
\phi_p(x_i)=\phi_p(x_i+hf(x_i,a^*(x_i)) & x_i \in G\cap \Omega\,,
\end{array}
\end{equation}
which is a semi-Lagrangian discretization of the advection equation $(\ref{advection})$. 
Note that, due to the interpolation and differently from the continuous case, 
scheme $(\ref{discrete-advection})$ spreads the sharp values $\{0,1\}$ of the initial guess, 
producing a final solution $\phi_p$ valued on the whole interval $[0,1]$. 
This makes the patches \emph{fuzzy} sub-sets of $\Omega$. 
Hence, an additional \emph{thresholding} procedure is needed to define them correctly. Indeed, for a fixed $\tau_P\in(0,1)$, we set
$$
\Omega_p=\left\{ x_i\in G\cap \Omega\,:\, \phi_p(x_i)\geq \tau_P\right\}\,.
$$
In practice, we loose some information, projecting $\phi_p$ on the space of characteristic functions. 
Nevertheless, this allows to choose a poor tolerance $\varepsilon_P$, accelerating the 
convergence of the scheme. Moreover, by tuning the threshold parameter $\tau_P$, we can achieve a desired level of overlapping 
between the patches, including sharp interfaces. We recall that, by condition $(\ref{upper-bound})$, 
the stencil of the semi-Lagrangian scheme $(\ref{SLscheme3})$ consists of first neighboring nodes. In addition, 
in this more general setting, we can no longer impose state constraint conditions at the boundaries of the patches. 
Then, some overlap is required, if we work with distributed memory architectures. On the other hand, sharp 
interfaces are still possible, if we prefer a shared memory architecture, as the one employed for the numerical experiments 
presented in the next section.
\\

\noindent Finally, the new patchy algorithm summarizes as follows:\\

Initialization:
\begin{itemize}
 \item Build two grids $G_c$ and $G$ such that $G_c<<G$.
 \item Fix tolerances $\varepsilon$, $\varepsilon_c$, $\varepsilon_P$, 
 the number of patches $N_P$ and the threshold $\tau_P$. 
 \item Build the initial guess $u_c^0$ on $G_c$ equal to the exit cost $g$ 
 on $G\cap\partial\Omega$ and $+\infty$ (i.e. a very big value) otherwise.
 \end{itemize}
 
 Pre-computation:
 \begin{itemize}
 \item Starting from $u_c^0$, compute $u_c$ on $G_c$, 
 iterating the scheme $(\ref{SLscheme3})$ with  $\sigma\equiv 0$ until convergence (up to $\varepsilon_c$).
 \item Build $\hat u$ on $G$ by interpolation of $u_c$.
 \item Compute $\hat a^*$ on $G$ using $\hat u$ in synthesis procedure $(\ref{discrete-synthesis})$.
 \item Divide the boundary nodes of $G$ in $N_P$ sub-sets.
 \item For $p=1,...,N_P$, use $\hat a^*$ to compute on $G$ the patch $\Omega_p$ (with threshold $\tau_P$), iterating the scheme 
 $(\ref{discrete-advection})$ until convergence (up to $\varepsilon_P$).
 \item For $p=1,...,N_P$, sort the nodes of $\Omega_p$ according to the increasing values of $\hat u$ (to exploit causality).
\end{itemize}
 
 Computation:
 \begin{itemize}
 \item For $p=1,...,N_P$, starting from $\hat u$, compute $u_p$ on $\Omega_p$, iterating the scheme 
 $(\ref{SLscheme3})$ until convergence (up to $\varepsilon$). At each iteration, update $u_p$ at the nodes where $\Omega_p$ 
 intersects other patches, using the transmission condition 
 $$
 \qquad u_p(x_i)=\min\{u_p(x_i),u_q(x_i)\}\mbox{ for all } q\neq p \mbox{ such that } \emptyset\neq\Omega_p\cap\Omega_q\ni x_i
 $$ 
 \item Build the solution $u$ on $G$, merging all the $N_P$ solutions $u_p$ on $\Omega_p$.
\end{itemize}
Note that all the steps of the method can be parallelized. In particular, the solution $u_c$ on the coarse grid $G_c$ can be 
computed by means of a standard domain decomposition technique. 
In the following section we will compare this dynamic domain decomposition to 
the classical static domain decomposition. The interested reader can find in \cite{QV08} a good introduction to this topic and in \cite{FLS94,CFLS94} some 
static domain decomposition methods for first order Hamilton-Jacobi equations.

We conclude this section by remarking that, in the case of second order semi-linear equations, 
the issue of balancing the size of the patches can no longer be addressed 
via the multi-level approach discussed in Section \ref{PD1} for first order equations. 
Indeed, due to the presence of the diffusion, we are not guaranteed that the values of the solution at 
the boundary of the current level decomposition are correct. Some information could flow back in the future from 
patches that have not yet been built. In principle, we can continue the construction of the decomposition, 
postponing the computation of the solution. To this end we need more processors, exactly $N_P N_L$, 
where $N_L$ is the number of levels, which is clearly bounded but a priori unknown. 
In the case of only $N_P$ available processors, an alternative could be to find an iterative method to solve 
the following optimization problem: 
build an initial subdivision of the boundary $\partial \Omega$ such that the corresponding patches have about the 
same sizes. An interesting question is to understand if the additional 
computational cost of this optimization procedure is compensated or not by the balancing in size of the patches. 
This method is at present under development.

\section{Numerical experiments}\label{experiments}
In this section, we present some numerical experiments in dimension two, 
performed on a server Supermicro 8045C-3RB using 1 CPU Intel Xeon Quad-Core E7330 2.4 GHz with 32 GB RAM, running under the Linux
Gentoo operating system. The aim is to emphasize the features of the 
proposed method. In particular, we first show that the modified semi-Lagrangian scheme \eqref{SLscheme3} 
allows for a substantial reduction of the number of iterations needed to reach convergence, compared to the 
original scheme \eqref{SLscheme2}. 
Next we show that it is able to compute the solution to general nonlinear equations with very degenerate diffusion. 
Finally, we combine the scheme with the upwind diffusion ball condition \eqref{hyperbolicity} 
and the fast-marching-like sorting of the grid nodes, in order to get an extra reduction of the computational time.  
This will confirm the effectiveness of the extension of the patchy method to the more general setting of second 
order Hamilton-Jacobi equations. 

In all the following tests, we set the domain $\Omega=[-1,1]^2$ and the control set $A=B_1$ 
(the unit ball centered in the origin), which is  discretized by means of $16$ points. 
Moreover, if not differently specified, we take the boundary datum $g(x)\equiv 0$, 
the running cost $l(x,a)\equiv 1$ and 
the diffusion $\sigma(x,a)\equiv \sqrt{2\varepsilon}{I}_2$, where $\varepsilon\ge0$ 
and ${I}_2$ denotes the $2\times2$ identity matrix. The corresponding second order operator in equation \eqref{HJB2} 
is just $-\varepsilon\Delta u$, i.e. the laplacian with diffusion coefficient $\varepsilon$. 
Finally, we denote by $\chi_S$ the characteristic function of a generic subset $S$ of $\Omega$.

Now, let us give a list of typical problems, depending on the choice of the dynamics $f(x,a)$:

\begin{itemize}
 \item[\bf A)] $f(x,a)=b(x)$ for a given vector field $b:\Omega\to\mathbb{R}^2$. 
 The corresponding equation is a stationary advection equation 
 along $b$ with uniform source: 
 $$b(x)\cdot\nabla u(x)=1\,.$$
\item[\bf B)] $f(x,a)=c(x)a$ for a given positive speed function $c:\Omega\to\mathbb{R}$. 
In this case it is easy to see that the $\max$ in \eqref{HJB1} is achieved for 
$a=-\nabla u(x)/|\nabla u(x)|$. The corresponding equation is the eikonal equation, whose solution $u(x)$
represents the minimum time to reach the boundary $\partial \Omega$ starting from $x\in\Omega$ and traveling 
at speed $c$: 
$$c(x)|\nabla u(x)|=1\,.$$
For $c\equiv 1$ we recover the distance function from the boundary $\partial\Omega$.\\
\item[\bf C)] $f(x,a)=\frac{1}{1+|x|^2}\left(R_\theta \frac{x}{|x|}+\frac{\eta}{2} a\right)$, 
where $\eta\in\{0,1\}$ is a fixed switch parameter to activate/deactivate the control and 
$R_\theta$ is a fixed counter-clock-wise rotation with $\theta<\frac{\pi}{2}$. The corresponding equation is 
a controlled advection equation with uniform source, also termed the Zermelo navigation problem. 
This is an example hard to compute, due to the strong anisotropy introduced by the control, namely the dependency of the speed 
$|f(x,a)|$ on the direction $a$. See below for further details.
\end{itemize}
Whenever $\varepsilon\neq 0$ we use to add the term \emph{diffusion} to declare the presence of the laplacian, as for the 
following eikonal-diffusion equation
\begin{equation}\label{eikdiff}
\left\{ \begin{array}{ll}
         -\varepsilon\Delta u(x)+c(x)|\nabla u(x)|=1 & x\in \Omega\\
         u(x)=0 & x\in \partial \Omega
        \end{array}
\right.
\end{equation}
\subsection{Self-dependency removal}
\begin{figure}[htp]
 \begin{center}
 \begin{tabular}{cc}
  \includegraphics[width=0.3\textwidth]{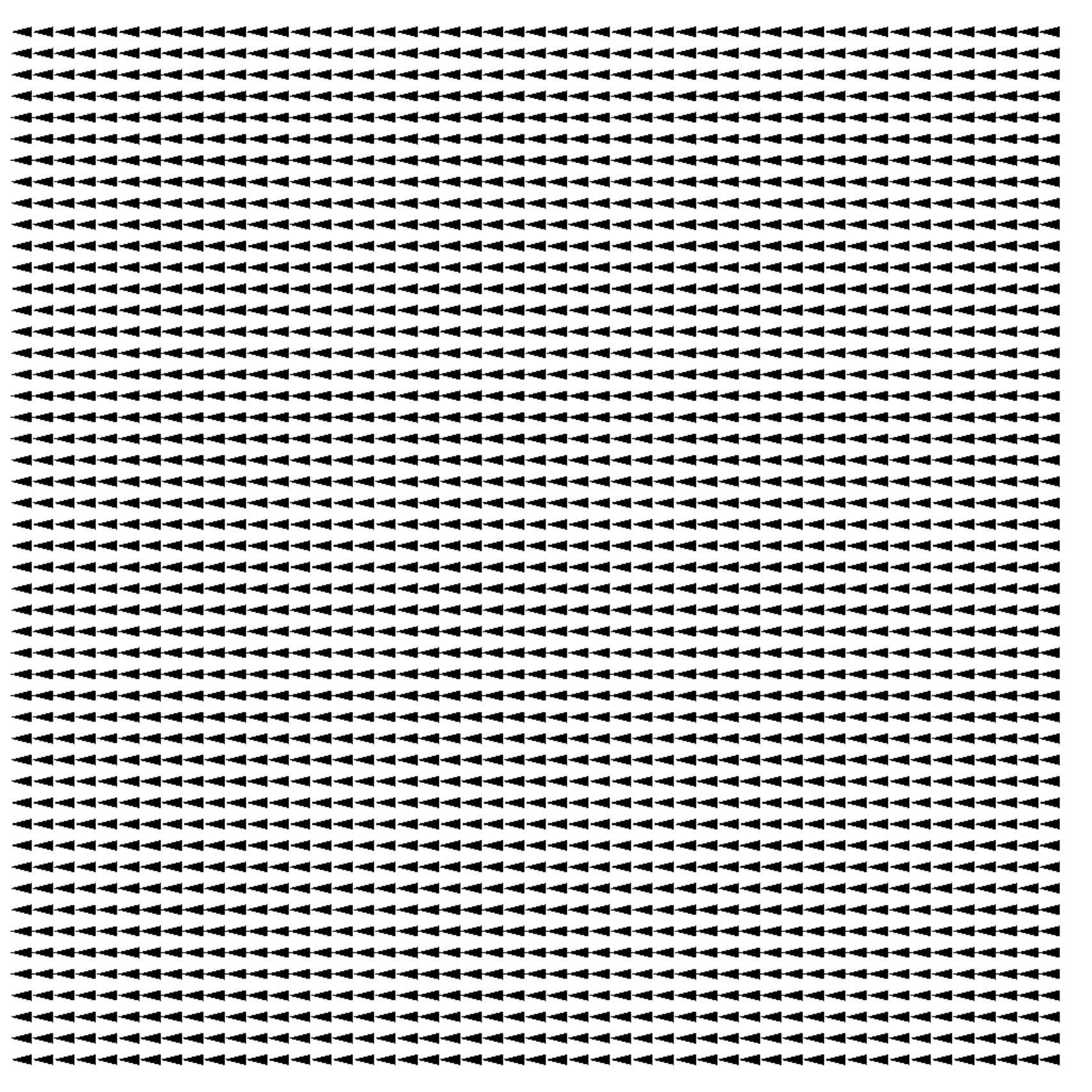} 
  &
  \includegraphics[width=0.3\textwidth]{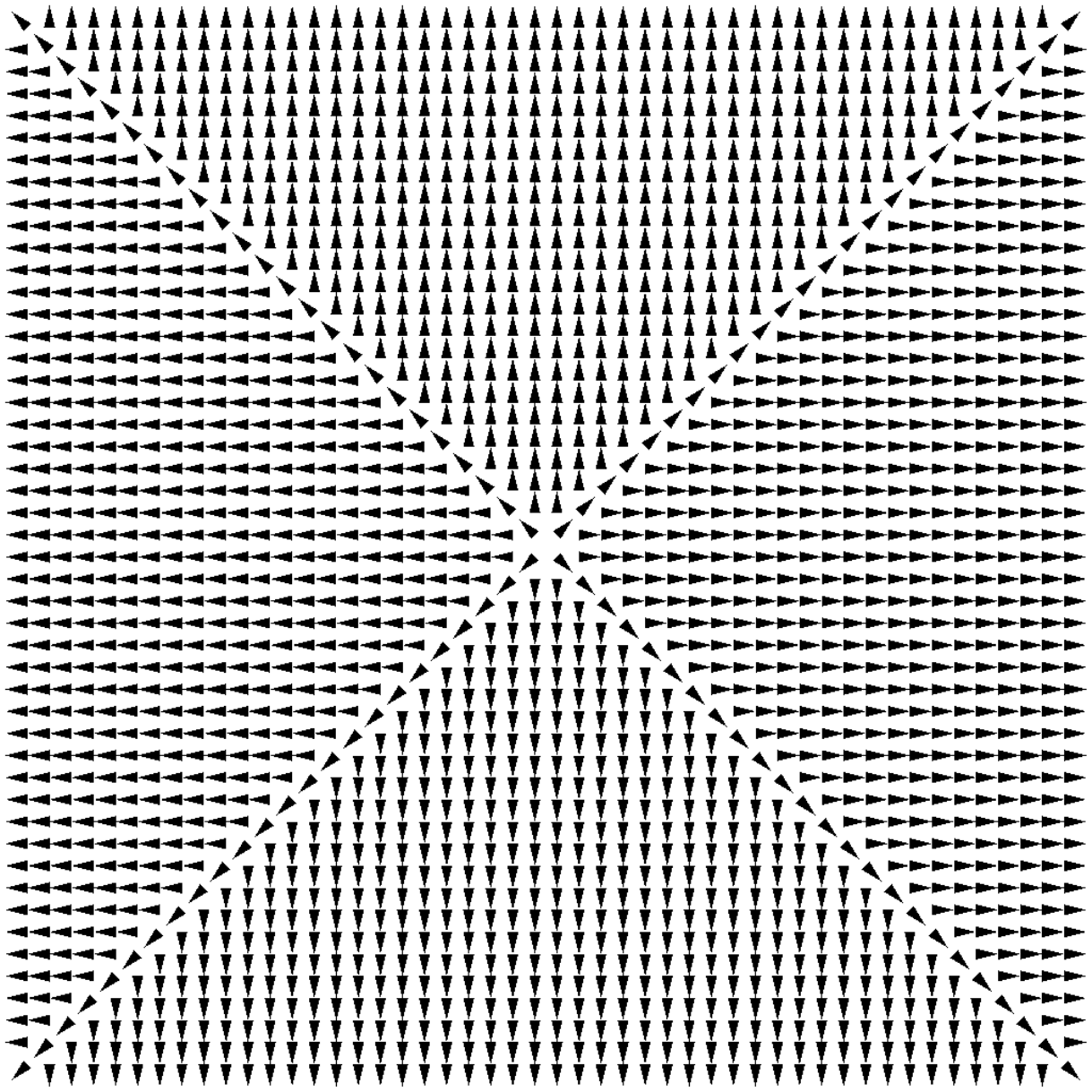} \\
  (a) & (b) 
 \end{tabular}
  \caption{Optimal vector fields for the stationary advection equation (a) and the eikonal equation (b).}\label{causality}
   \end{center}
 \end{figure}

Here we compare our modified semi-Lagrangian scheme 
\eqref{SLscheme3} and the original one \eqref{SLscheme2} without any removal of the self-dependency. 
In Table \ref{table0} we report the number of iterations needed to converge for both schemes 
in the cases listed above. All the solutions are computed on a $50\times 50$ grid. 
\begin{table}[!h]
 \centering
  \footnotesize
   \caption{The modified SL scheme \eqref{SLscheme3} vs the original one \eqref{SLscheme2}. 
  Number of iterations to reach convergence
  for different test dynamics and diffusion coefficients.}\label{table0} 
  \begin{tabular}{|c|c|c|c|c|}
    \hline
    Equation & Dynamics & $\varepsilon$ & SL \eqref{SLscheme2} & SL \eqref{SLscheme3}\\
    \hline
    A & $b(x)\equiv (1,0)$ & $0$ & 176 & 2\\  
     \hline
    A & $b(x)\equiv (1,0)$ & $1e\mbox{-}2$ & 458 & 102\\  
     \hline
    B & $c(x)\equiv 1$ & $0$ & 145 & 26\\  
     \hline
    B & $c(x)\equiv 1$ & $1e\mbox{-}2$ & 353 & 98\\  
     \hline
    B & $c(x)=1+\chi_{\{x_1\ge 0\}}(x)$ &
  $0$ & 276 & 33\\  
   \hline
 B & $c(x)=1+\chi_{\{x_1\ge 0\}}(x)$ &
  $1e\mbox{-}2$ & 305 & 78\\  
   \hline
 C & $\eta=1$ & $0$ & 716 & 116\\ 
   \hline
 C & $\eta=1$ & $1e\mbox{-}2$ & 961 & 325\\ 
\hline
 \end{tabular}
\end{table}\\
For a pure advection dynamics (A) with vector field $b=(1,0)$ and $\varepsilon=0$, 
the causality property is explicit, in the sense that, starting from the left boundary, information propagates 
to the right in the whole domain. Since, by default, we process the grid nodes from left to right and from top to 
bottom, then the modified SL scheme converges in this case in just one iteration. Each node is computed (backward in time) using 
only nodes on its left, hence it is correct by induction (see Figure \ref{causality}a). 
Note that the additional iteration reported in the 
table is the one needed to check the convergence. On the contrary, the original scheme employs 
the value of a node (which is set at the beginning to a big value as initial guess) 
to compute itself. 
It follows that the correct value carried by the left neighbor is changed by the interpolation and several additional 
fixed-point iterations are needed to fix it. 
The case of the eikonal equation (B) with uniform speed $c=1$ and $\varepsilon=0$ is similar, as shown in Figure \ref{causality}b. 
It turns out that characteristics are straight lines, moving from the boundary until they intercept the diagonals of the square. 
Using as before the default order for processing the grid nodes (left to right, top to bottom), 
it is easy to see that, at the end of the first iteration, only the nodes above the diagonal $x_1=x_2$ 
will contain correct values of the solution. For the remaining part of the domain, 
only the first neighboring nodes of the boundary are correctly updated. The same holds in the following iterations, 
as long as information reaches the \emph{middle} of the $50\times 50$ grid. This gives precisely $25$ 
iterations (plus one to check the convergence) for the modified SL scheme and much more (145) for the original 
scheme. This behavior is confirmed also for more complicated dynamics, as for the eikonal equation (B) with a 
non homogeneous speed function and for the anisotropic dynamics (C). We finally remark that, in presence of diffusion $\varepsilon>0$, 
the causality property of hyperbolic equations is lost and the number of iterations 
necessarily increases. Nevertheless, we still get a relevant reduction of this number for the modified SL scheme.
From now on, only this scheme will be employed.
\subsection{Degenerate diffusion}
\begin{figure}[htp]
\begin{center}
 \includegraphics[width=0.3\textwidth]{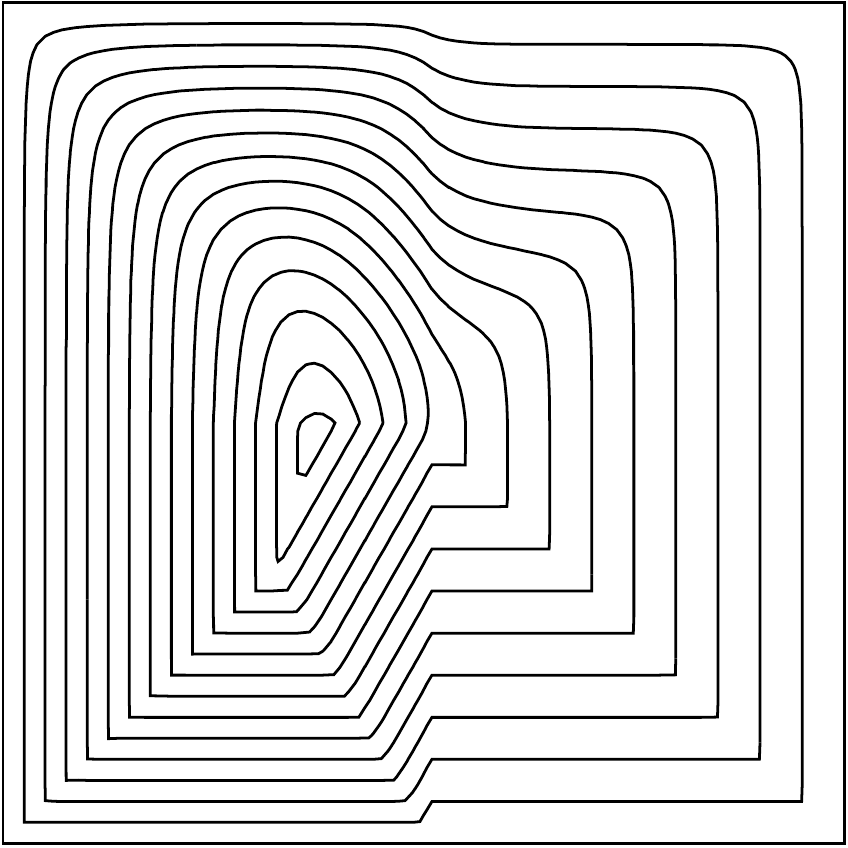}\\
 \caption{Level sets of the solution to a degenerate eikonal-diffusion equation with a non homogeneous speed: the diffusion is 
 present only in the upper half of the square, 
 affecting the smoothness of the solution.}\label{degeneracy}
  \end{center}
\end{figure}

The following tests aim at showing the built-in ability of the proposed semi-Lagrangian scheme to solve problems 
with degenerate diffusion. To this end, we allow the diffusion coefficient $\varepsilon$ to depend also on $x$, taking the form 
$\varepsilon(x)=0.1\chi_{\{x_2\ge 0\}}(x)$. 
With this choice, we consider again the eikonal-diffusion equation \eqref{eikdiff} with speed $c(x)=1+\chi_{\{x_1\ge 0\}}(x)$
and we note that it is \emph{elliptic} in $\Omega\cap{\{x_2\ge 0\}}$ and \emph{hyperbolic} otherwise. 
Accordingly, at the grid nodes belonging to $\Omega\cap{\{x_2< 0\}}$, we recover by construction 
the semi-Lagrangian scheme designed for first order hyperbolic equations. Figure \ref{degeneracy} 
shows the level sets of the corresponding solution. 
We can immediately distinguish the two different behaviors, by looking at the smoothness of the solution. 
In particular, the sharp corners corresponding to the shocks in the hyperbolic part of the domain, 
are completely smoothed out in $\Omega\cap{\{x_2\ge 0\}}$. 
 \\

\begin{figure}[htp]
\begin{center}
\begin{tabular}{ccc}
 \includegraphics[width=0.3\textwidth]{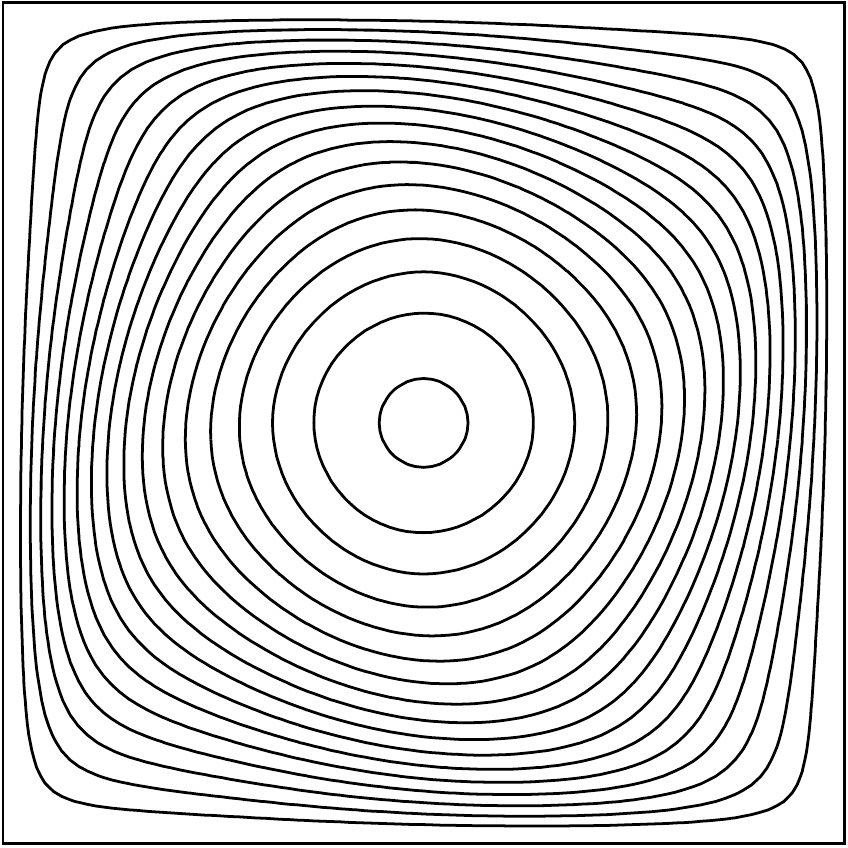} 
 &
 \includegraphics[width=0.3\textwidth]{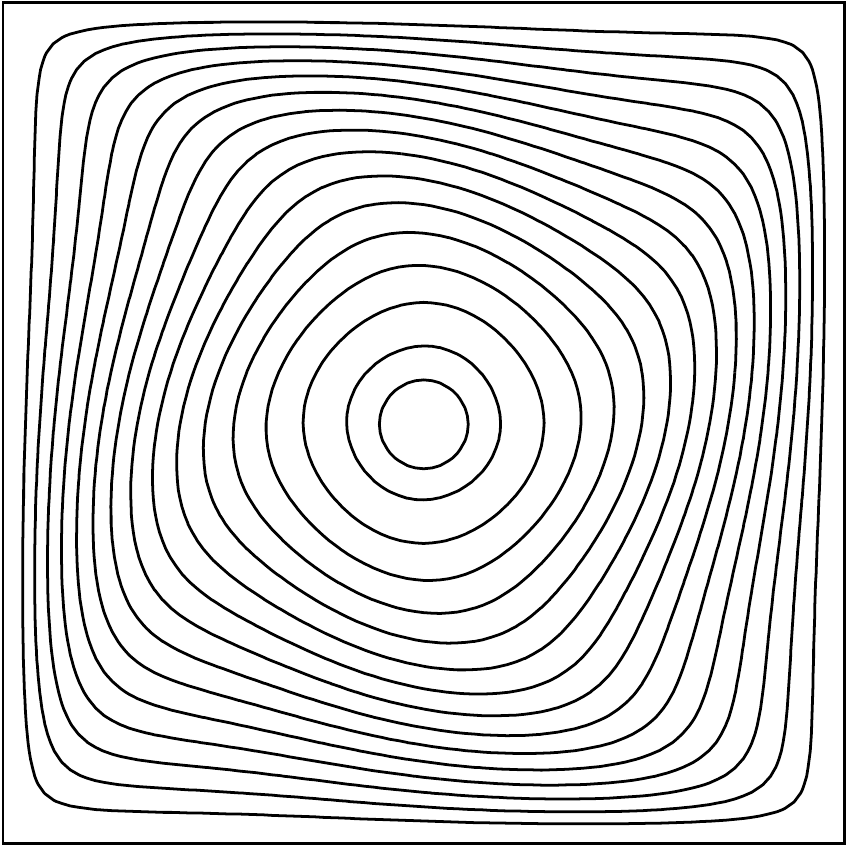} 
 &
 \includegraphics[width=0.3\textwidth]{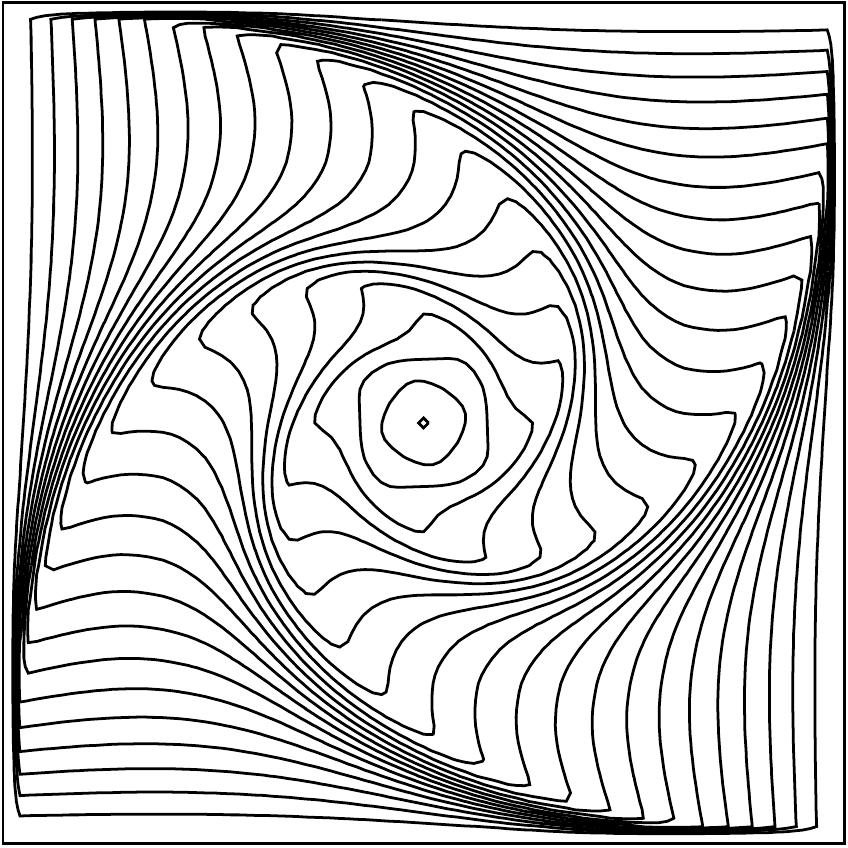}\\
 
 \includegraphics[width=0.3\textwidth]{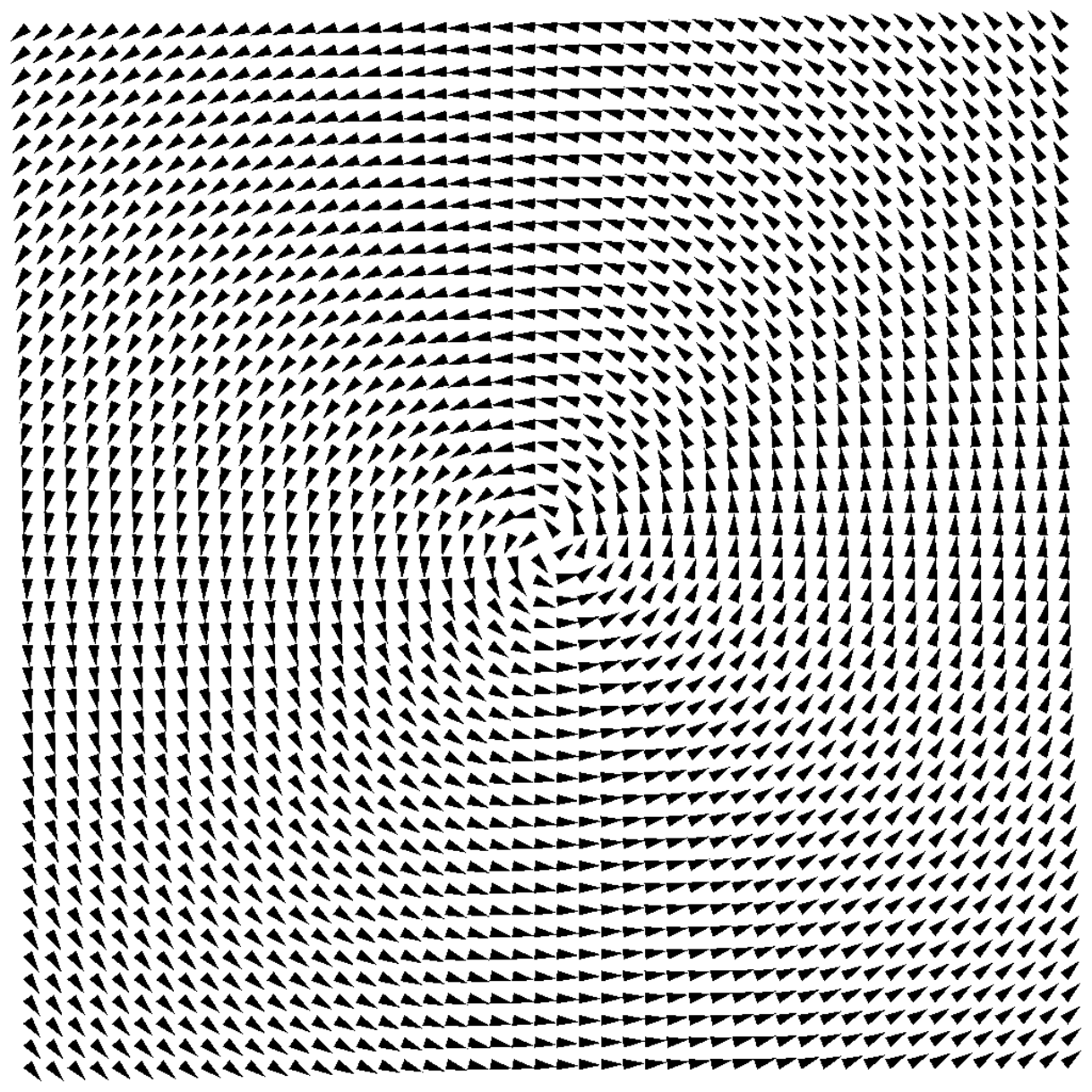} 
 &
 \includegraphics[width=0.3\textwidth]{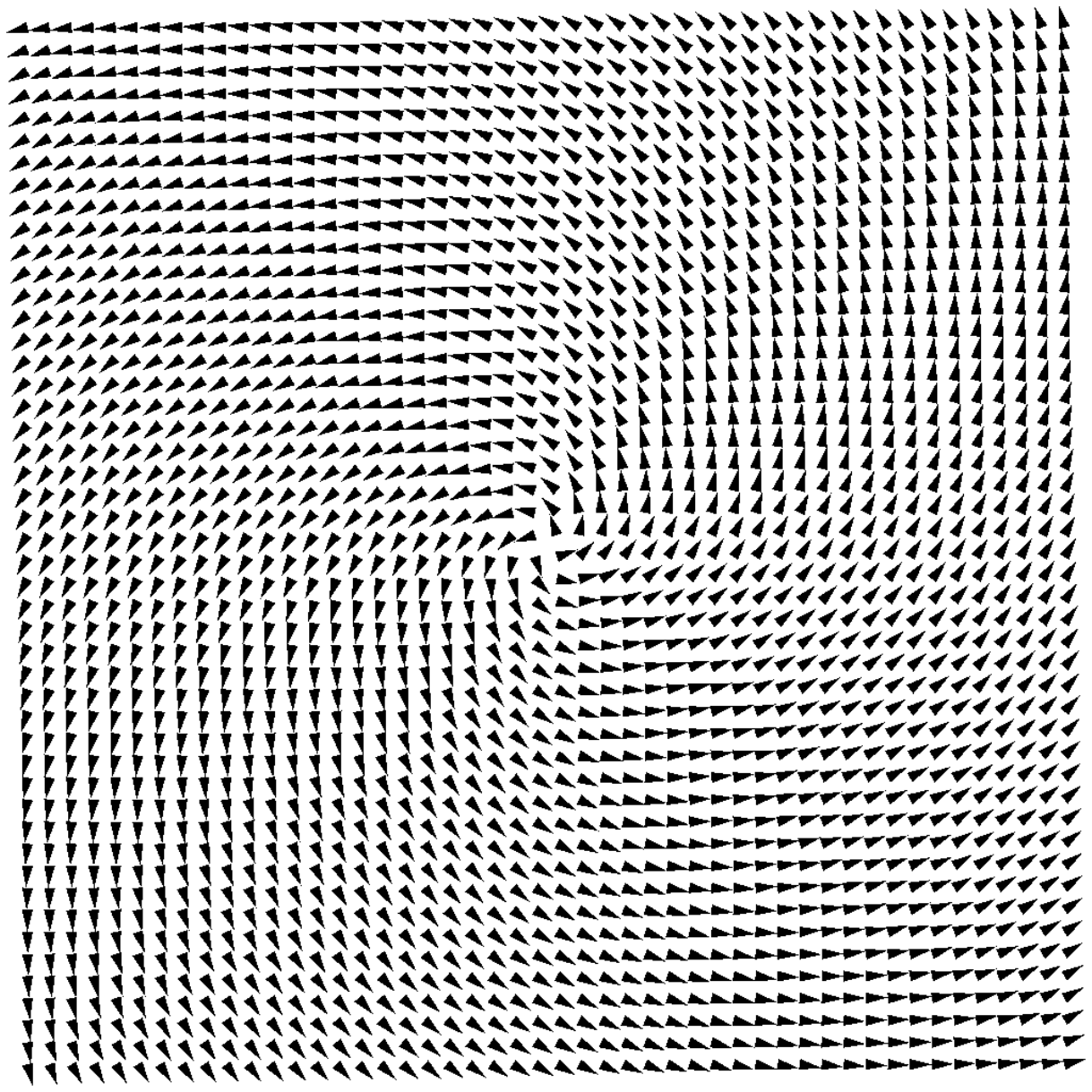} 
 &
 \includegraphics[width=0.3\textwidth]{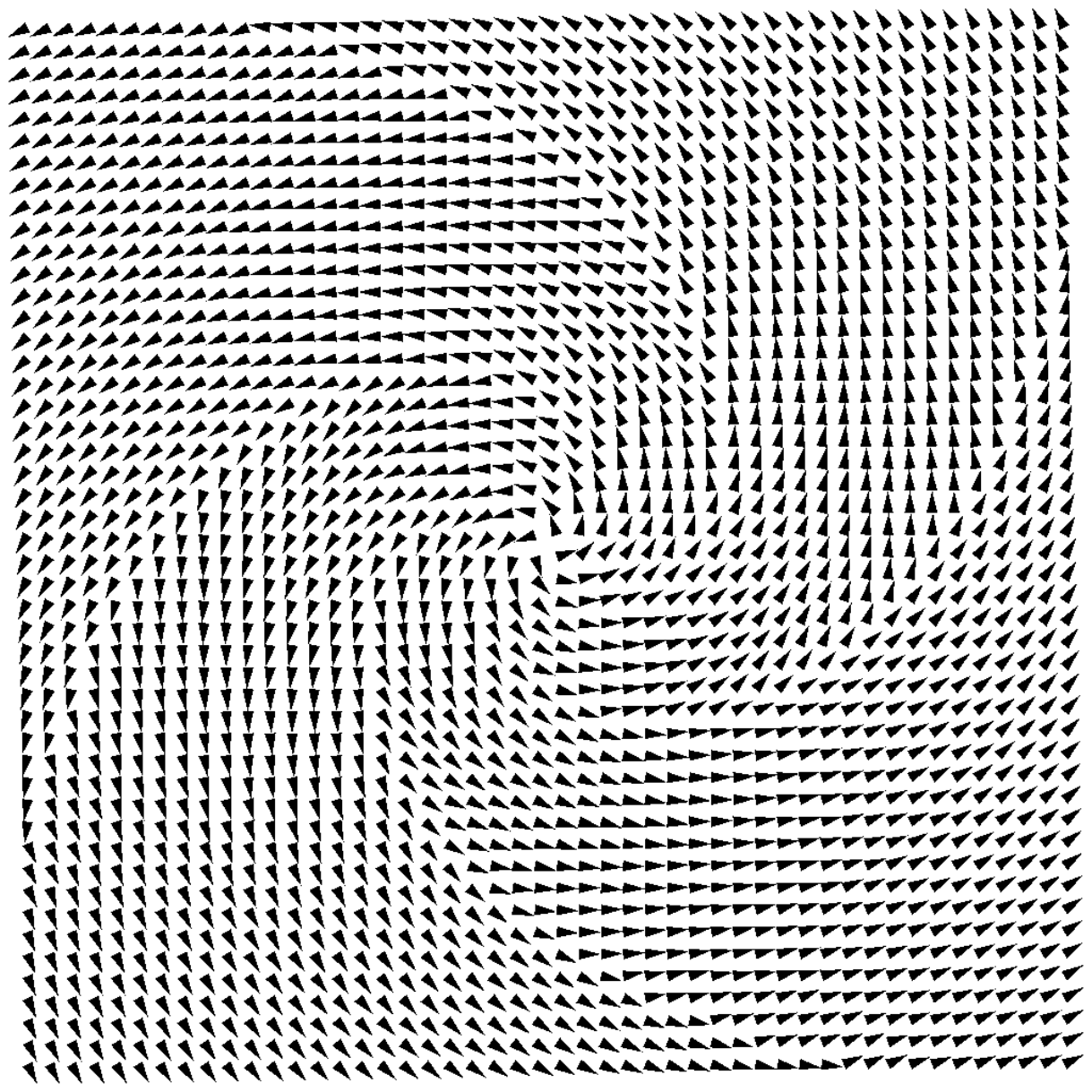}\\
 (a) & (b) & (c)
\end{tabular}
 \caption{The Zermelo navigation problem, level sets of the solution and optimal dynamics: 
 (a) diffusion without control ($\varepsilon=0.1$, $\eta=0$), (b) diffusion with control ($\varepsilon=0.1$, $\eta=1$), 
 (c) control without diffusion ($\varepsilon=0$, $\eta=1$).}\label{degeneracy2}
  \end{center}
\end{figure}

Now we consider the more suggestive and complicated example of the Zermelo navigation problem, whose dynamics is 
rewritten here for the reader's convenience:
$$f(x,a)=\frac{1}{1+|x|^2}\left(R_\theta \frac{x}{|x|}+\frac{\eta}{2} a\right)\,.$$
This problem consists in reaching the boundary $\partial\Omega$ in minimum time, starting from a point in $\Omega$ 
and using the control $a$. 
The first term in the dynamics represents a whirling drift from the origin to the boundary, 
whereas the fixed parameter $\eta\in\{0,1\}$ allows to switch on/off the control. 
This introduces in the problem both inhomogeneity and strong anisotropy.
It is easy to see that, for a fixed $\theta<\frac{\pi}{2}$, the counter-clock-wise rotation $R_\theta$ guarantees 
the reachability of the boundary for each starting point in $\Omega$, even for $\eta=0$, but we want to play with the control $a$ 
to minimize the time of arrival. In Figure \ref{degeneracy2} we compare 
the level sets of the solutions and the corresponding optimal dynamics for three different cases.
The first case (see Figure \ref{degeneracy2}a) represents a pure advection-diffusion equation, 
with uniform diffusion $\varepsilon=0.1$ on the whole $\Omega$ and no control, i.e. $\eta=0$. In the second case 
(see Figure \ref{degeneracy2}b) also the control is active, i.e. $\eta=1$, and we 
can clearly see how it produces a resistance to the drift that allows to reach the boundary in a smaller time. 
This is much more evident in the third case 
(see Figure \ref{degeneracy2}c) where the control is still active and the diffusion is switched off, i.e. $\varepsilon=0$.
\subsection{General nonlinear problems}
\begin{figure}[htp]
\begin{center}
\begin{tabular}{ccc}
 \includegraphics[width=0.3\textwidth]{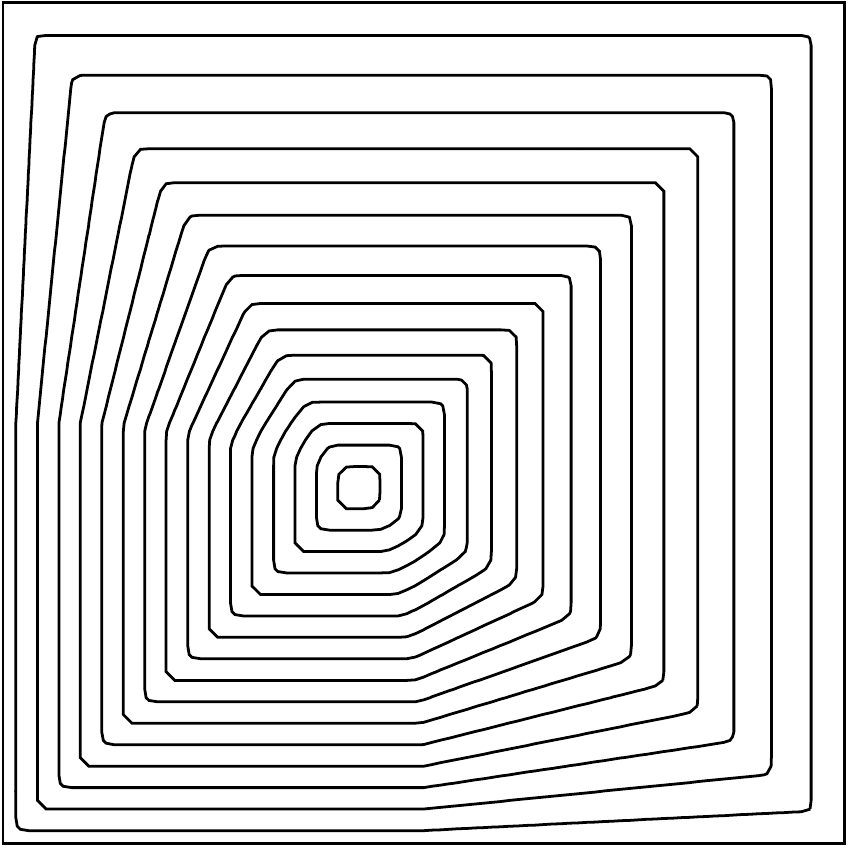} 
 &
\includegraphics[width=0.3\textwidth]{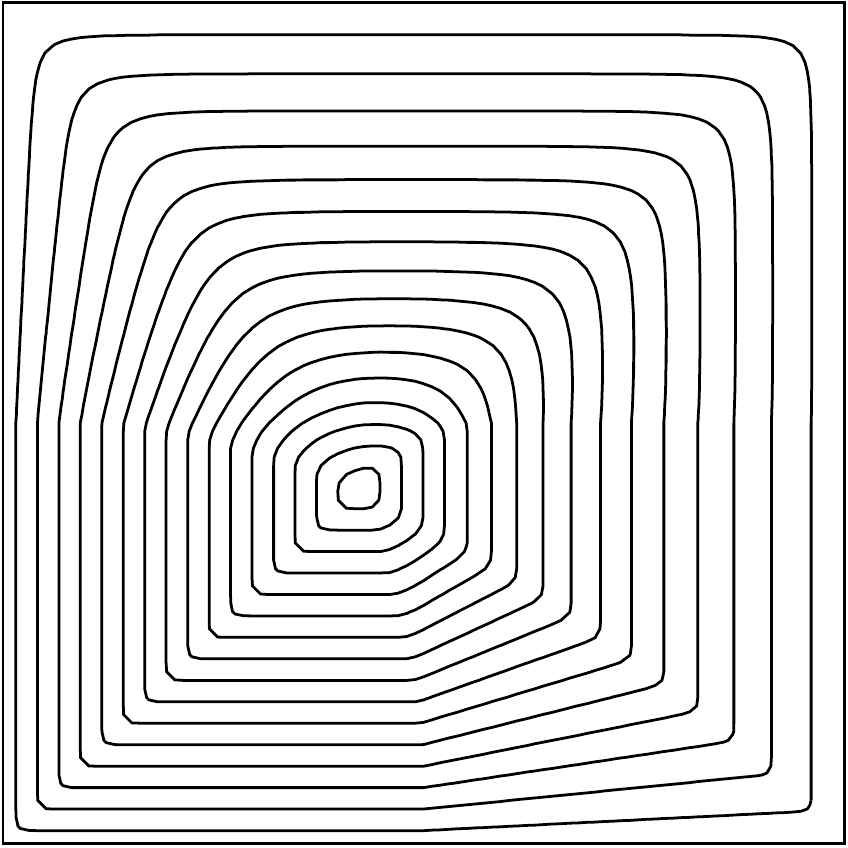} 
 &
 \includegraphics[width=0.3\textwidth]{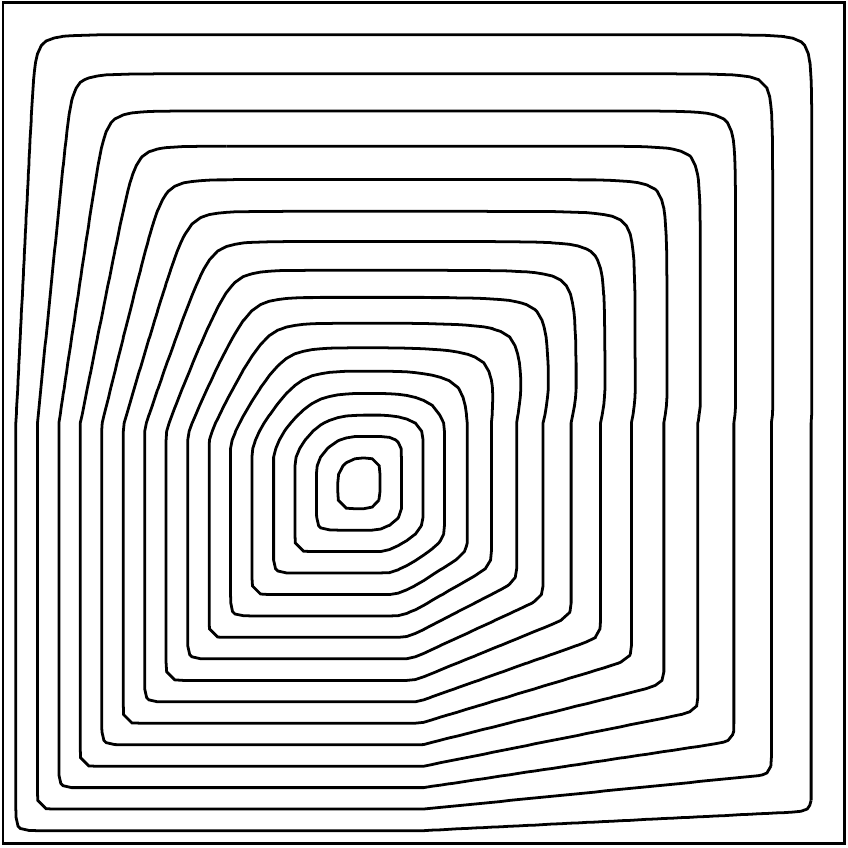} 
 \\
 $(l_1,\sigma_1)$ & $(l_1,\sigma_2)$ & $(l_1,\sigma_3)$
 \\
 \includegraphics[width=0.3\textwidth]{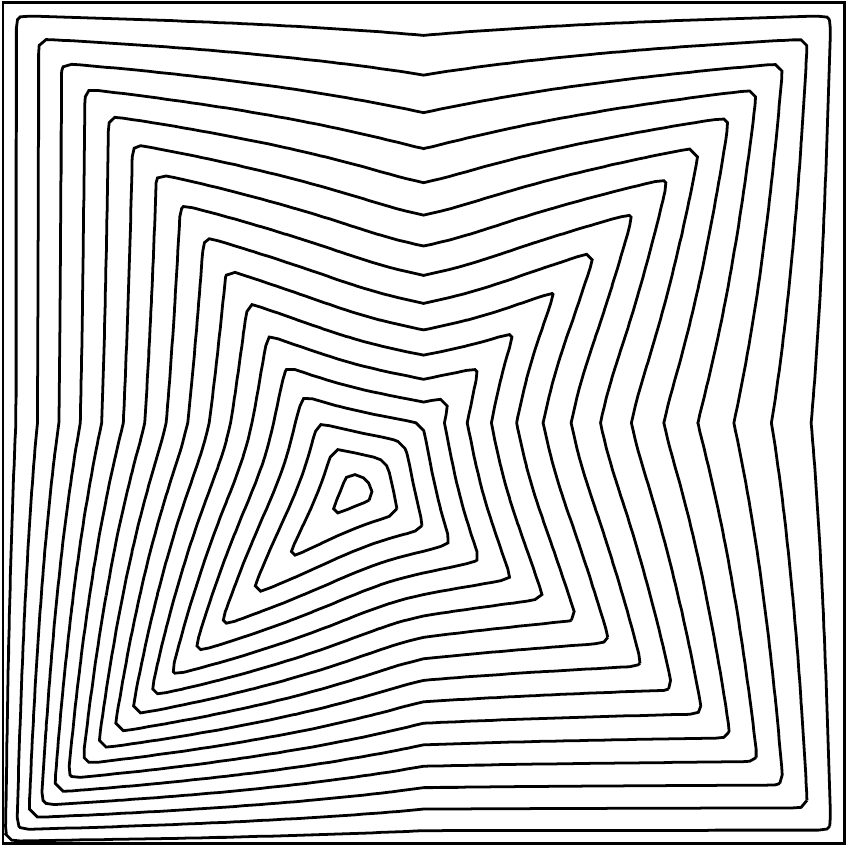} 
 &
\includegraphics[width=0.3\textwidth]{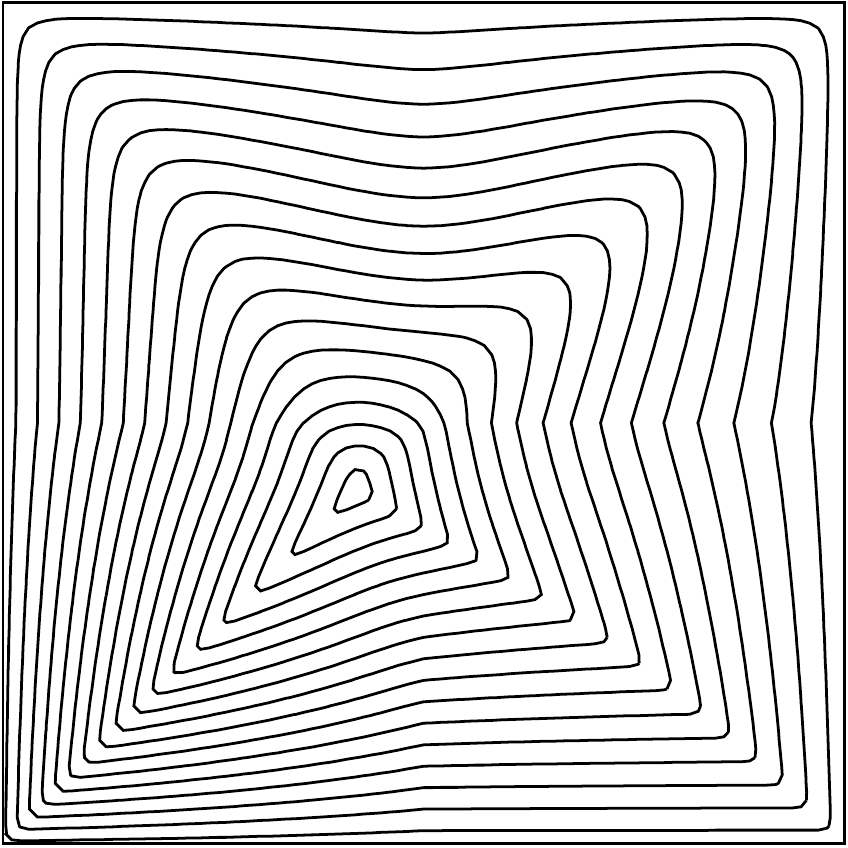} 
 &
 \includegraphics[width=0.3\textwidth]{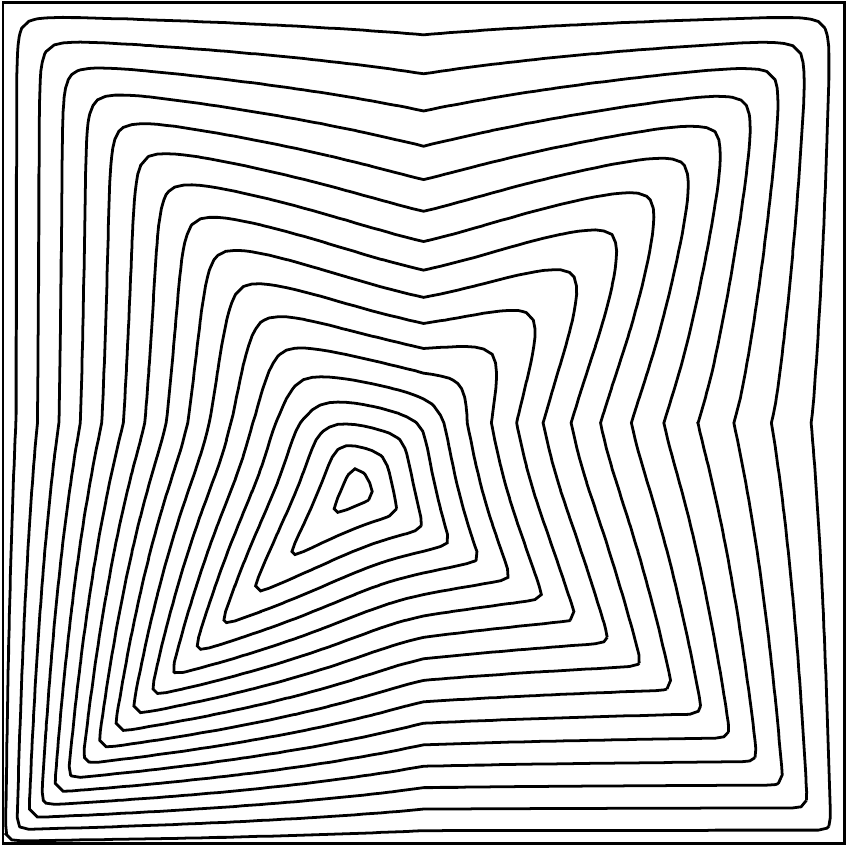} 
 \\
 $(l_2,\sigma_1)$ & $(l_2,\sigma_2)$ & $(l_2,\sigma_3)$
 \\
 \includegraphics[width=0.3\textwidth]{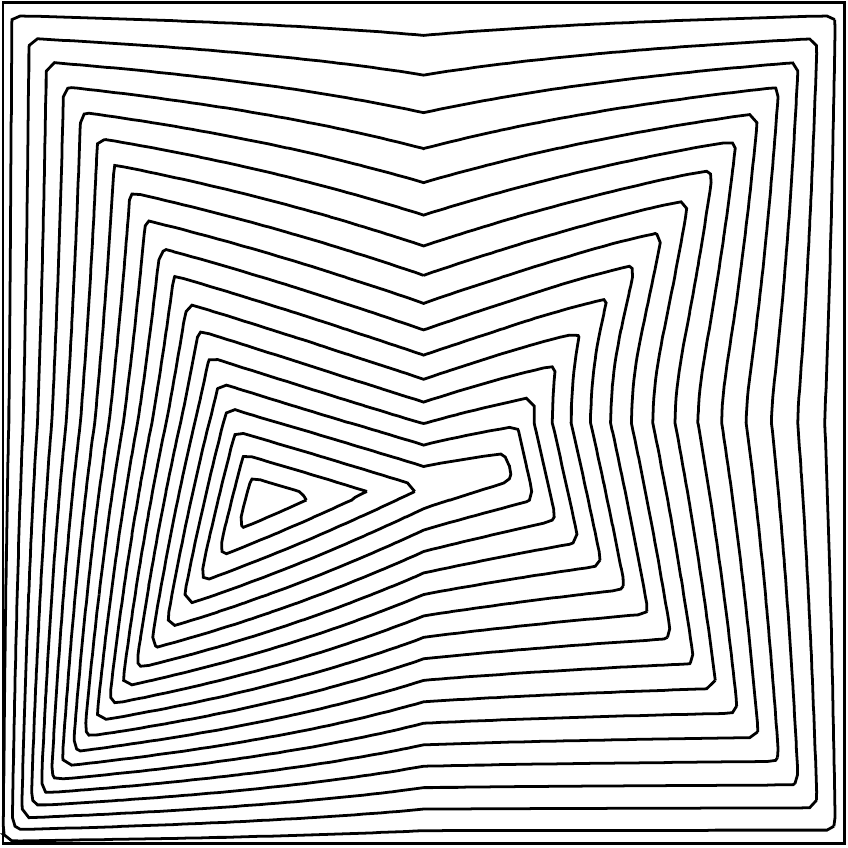} 
 &
\includegraphics[width=0.3\textwidth]{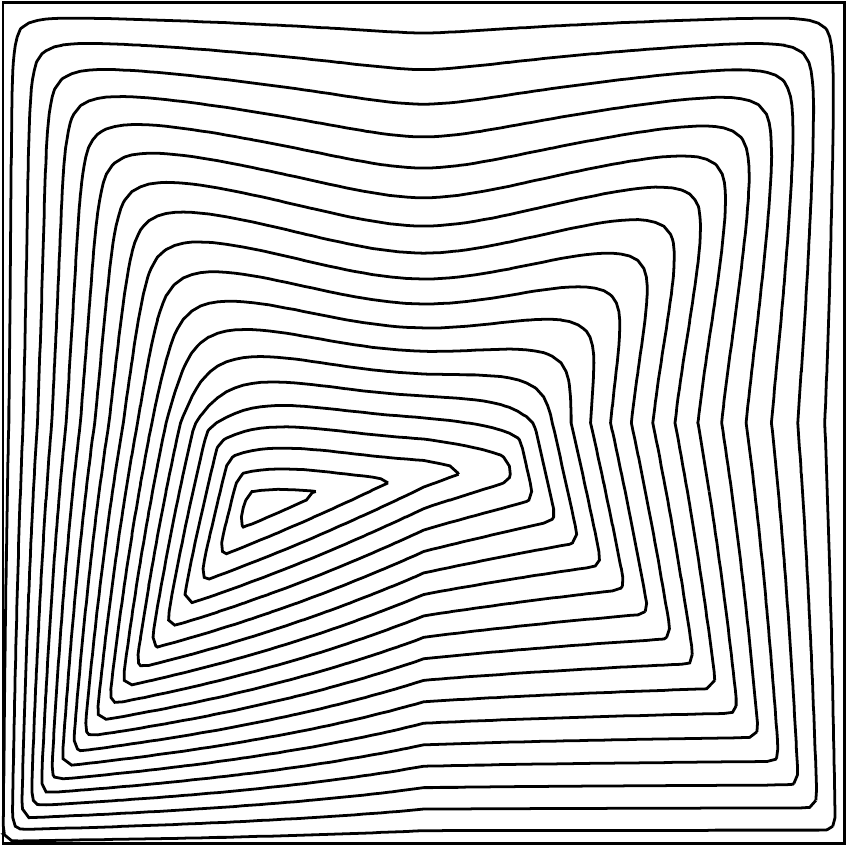} 
 &
 \includegraphics[width=0.3\textwidth]{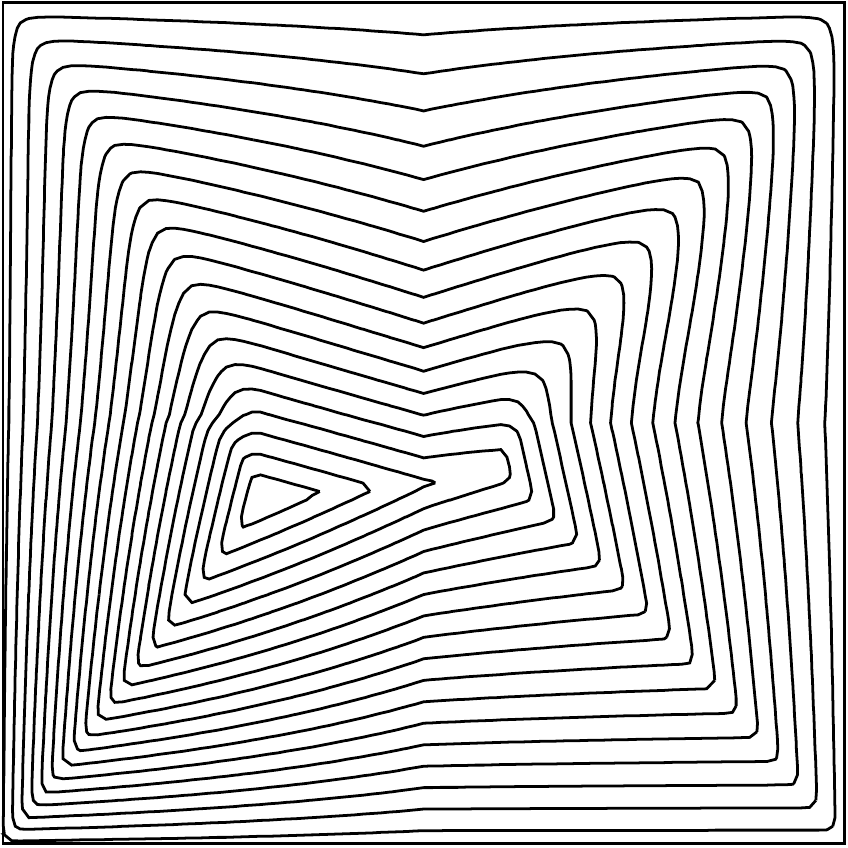} 
 \\
 $(l_3,\sigma_1)$ & $(l_3,\sigma_2)$ & $(l_3,\sigma_3)$
 \end{tabular}
 \caption{Level sets of the solutions corresponding to the dynamics \eqref{nonomog} and 
 to the pairs $(l_i,\sigma_j)$ ($i,j=1,2,3$) of running costs \eqref{l123} and diffusion terms \eqref{s123}.}\label{generalsigmacost}
  \end{center}
\end{figure}

This test shows the ability of the proposed semi-Lagrangian scheme to solve nonlinear problems with very general 
diffusion terms $\sigma$ and running costs $l$, possibly depending on $x\in\Omega$ 
and $a=(a_1,a_2)\in B_1$. To this end, we consider 
a non homogeneous dynamics of the form 
\begin{equation}\label{nonomog}
f(x,a)=c(x)a\quad\mbox{with}\quad c(x)=1+\max\{x_2,\max\{x_1,0\}\}\,.
\end{equation}
Moreover, we define 
$d_\varepsilon(x)=\sqrt{2\varepsilon}\chi_{\{x_2\ge 0\}}(x)$ and we consider three different diffusion terms:
\begin{equation}\label{s123}
\sigma_1\equiv0\,,\qquad \sigma_2(x)=d_\varepsilon(x)I_2\,,\qquad \sigma_3(x,a)=d_\varepsilon(x)a\,.
\end{equation}
Note that $\sigma_2$ corresponds to the usual uniform two dimensional Wiener process, whereas $\sigma_3$ defines a one dimensional 
Wiener process in the direction of the control $a$. Finally, we consider three different running costs:
\begin{equation}\label{l123}
l_1\equiv1\,,\qquad l_2(x)=1+|x_1 x_2|\,,\qquad l_3(x,a)=1+|x_1 x_2|+\left|\frac{a_1}{2+a_2}\right|\,.
\end{equation}
In Figure \ref{generalsigmacost} we show the level sets of the solutions corresponding to all the nine possible 
pairs $(l_i,\sigma_j)$ for $i,j=1,2,3$. We clearly see the effect of the diffusion, 
whenever applied in the upper half of the square. Moreover, we can also distinguish the different behavior between the 2D uniform diffusion and 
the 1D control driven diffusion. 
\subsection{Performance comparison}
\begin{figure}[htp]
\begin{center}
\begin{tabular}{ccc}
 \includegraphics[width=0.3\textwidth]{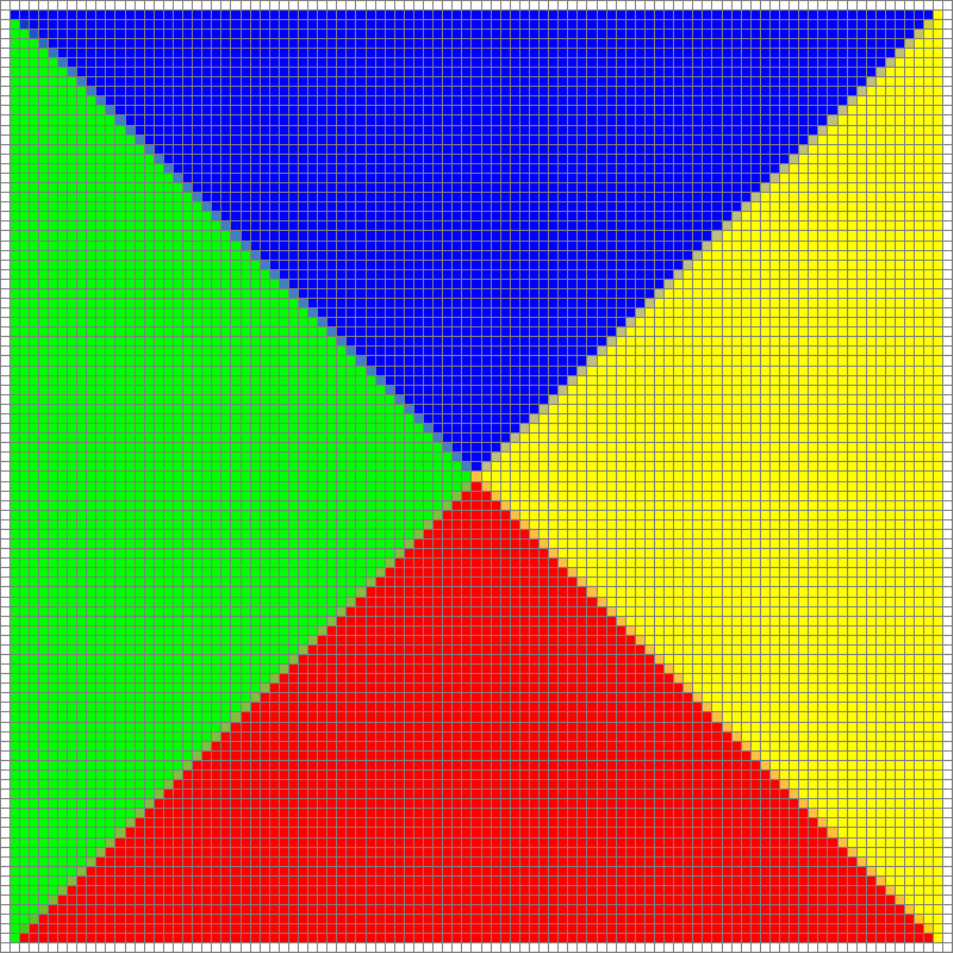} 
 &
\includegraphics[width=0.3\textwidth]{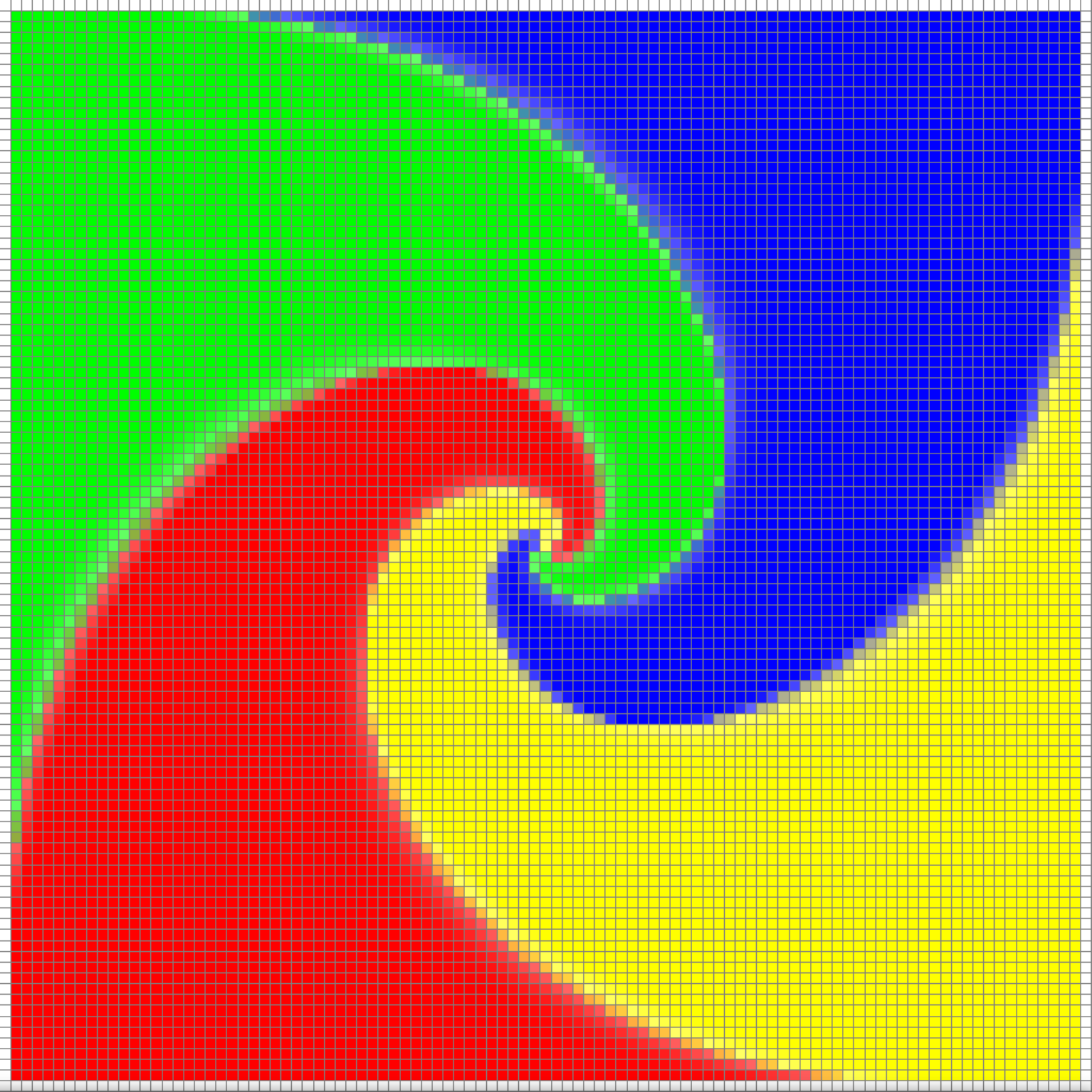} 
 &
 \includegraphics[width=0.3\textwidth]{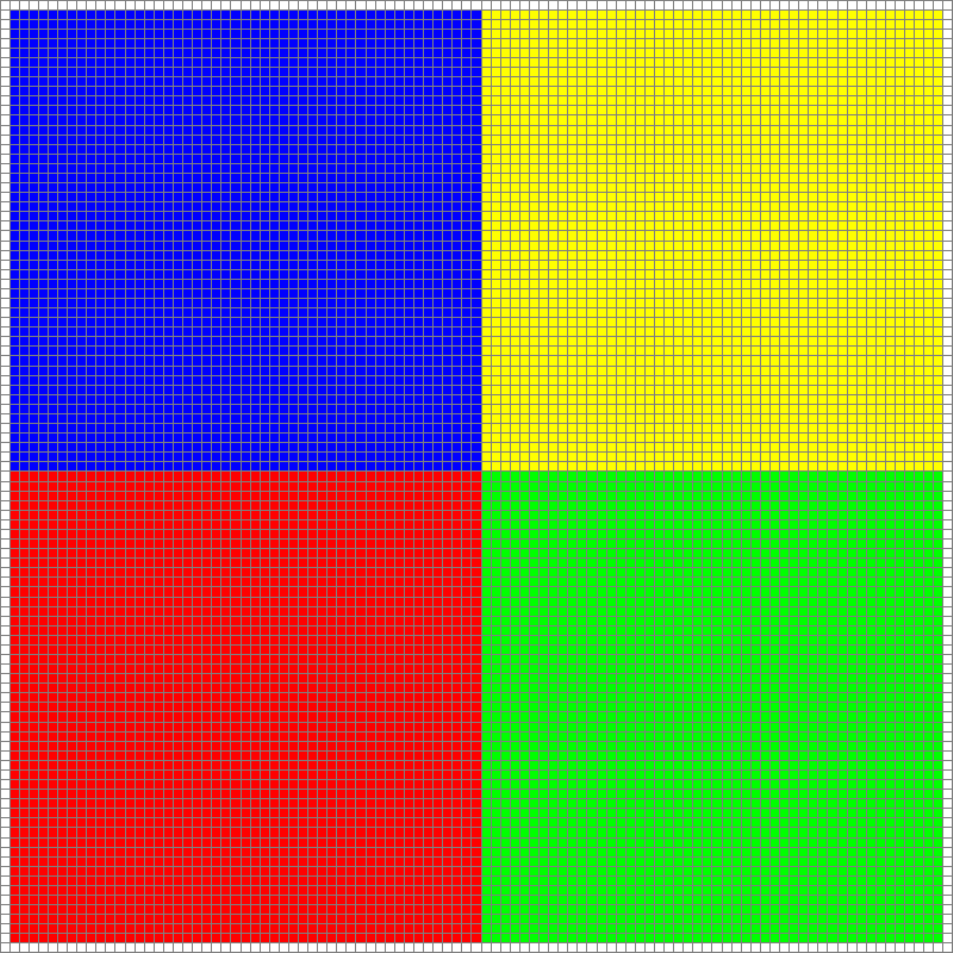} 
 \\
 (a) & (b) & (c)
\end{tabular}
 \caption{Domain decompositions in $4$ sub-domains: (a) PDD for the eikonal dynamics, 
 (b) PDD for the Zermelo dynamics, (c) a standard DD for both dynamics.}\label{decomposition}
  \end{center}
\end{figure}

Here we present a performance comparison between a standard Domain Decomposition method (DD) and the Patchy Domain 
Decomposition method (PDD). 
Here the aim is to convince the reader that, in order to obtain an additional speedup in the computation, 
the fundamental ingredients are the causality property and the upwind diffusion ball condition discussed in Section \ref{PD2}.
To this end, we consider the eikonal-diffusion equation with speed $c(x)\equiv 1$ and the Zermelo navigation problem above, 
both in the case with uniform diffusion $\varepsilon$ on the whole $\Omega$. We build the corresponding patchy domain 
decompositions starting from a subdivision of the boundary $\partial\Omega$ in four parts, namely the sides of the square. 
A natural choice to perform this computation in parallel is to employ four processors. 
Figure \ref{decomposition} shows, for a $100\times100$ grid, the resulting dynamic decompositions compared to an arbitrary and static one. 
Note that, in this case, the sub-domains of the three decompositions (the four triangles, the four spirals and the four squares respectively) 
have exactly the same size. 
 \\
We choose a starting coarse grid of $50\times 50$ nodes and different fine grids, up to $800\times 800$.
For each test, we sort the grid nodes of the patchy decompositions according to the increasing values of 
the pre-computed coarse solutions interpolated 
on the corresponding fine grid. As already remarked, for the eikonal equation ($\varepsilon=0$) this corresponds to the optimal order to exploit causality in 
a fast-marching fashion. 

In the following tables we report, for each fine grid and for different values of the diffusion coefficient 
$\varepsilon$, the CPU time in seconds and the number of iterations to reach convergence for the PDD method and a standard DD method. 
We remark that PDD's times also include the pre-computation step. 

\setlength{\tabcolsep}{1pt}
\begin{table}
 \centering
 \footnotesize
  \caption{The eikonal-diffusion equation: PDD vs DD in terms of CPU time in seconds (and number of iterations) 
  to reach convergence. Values in the first row of each cell concern PDD, the others concern DD. 
  Values in bold are obtained in the regime given by the upwind diffusion ball condition.}\label{table1} 
 \begin{tabular}{|c|c|cc|cc|cc|cc|c|}
 \cline{3-11}
 \multicolumn{2}{c|}{} & \multicolumn{2}{c|}{$100^2$} & \multicolumn{2}{c|}{$200^2$} & \multicolumn{2}{c|}{$400^2$} & \multicolumn{2}{c|}{$800^2$} & $N$\\ \cline{3-11}
 \multicolumn{2}{c|}{} &  \multicolumn{2}{c|}{0.02} & \multicolumn{2}{c|}{0.01} & \multicolumn{2}{c|}{0.005} & \multicolumn{2}{c|}{0.0025} & $\Delta x$ \\ \hline
 \multirow{18}*{$\varepsilon$}
 &\multirow{2}*{$10^{-9}$} & \bf 0.18&\bf(6)  & \bf0.62&\bf(6) & \bf3.39&\bf(6) &  \bf21.53&\bf(6) \\ 
 &                        & \bf0.57&\bf(52)  & \bf4.56&\bf(102) &  \bf35.88&\bf(202) & \bf297.98&\bf(402) \\ 
 \cline{2-10}

 &\multirow{2}*{$10^{-8}$} & \bf 0.22&\bf(9)  & \bf0.79&\bf(10) & \bf4.52&\bf(12) &  \bf27.48&\bf(14) \\ 
 &                        & \bf0.64&\bf(54)  & \bf4.61&\bf(104) &  \bf36.52&\bf(205) & \bf302.47&\bf(406) \\ 
 \cline{2-10}
  &\multirow{2}*{$10^{-7}$} & \bf 0.23&\bf(10)  & \bf0.89&\bf(12) & \bf4.92&\bf(14) &  \bf29.49&\bf(17) \\ 
 &                        & \bf0.68&\bf(55)  & \bf4.72&\bf(105) &  \bf37.53&\bf(206) & \bf304.59&\bf(408) \\ 
 \cline{2-10}

 &\multirow{2}*{$10^{-6}$} & \bf 0.26&\bf(12)  & \bf0.98&\bf(14) & \bf5.63&\bf(18) &  \bf34.76&\bf(23) \\ 
 &                        & \bf0.70&\bf(56)  & \bf4.77&\bf(107) &  \bf39.15&\bf(209) & \bf312.19&\bf(413) \\ 
 \cline{2-10}

 &\multirow{2}*{$10^{-5}$} &  \bf 0.28&\bf(14) & \bf1.19&\bf(18) & \bf6.94&\bf(25) & \bf42.82&\bf(35)\\
 &                        & \bf0.71&\bf(57) & \bf4.84&\bf(109) & \bf41.95&\bf(214) & \bf326.93&\bf(422)\\ 
 \cline{2-10} 
 
 &\multirow{2}*{$10^{-4}$} &  \bf 0.34&\bf(18) & \bf1.47&\bf(25) & \bf8.87&\bf(37) & \bf63.13&\bf(62)\\
 &                        & \bf0.79&\bf(60) & \bf5.08&\bf(114) & \bf44.18&\bf(223) & \bf349.71&\bf(445)\\ 
 \cline{2-10}
 &\multirow{2}*{$6.25\times 10^{-4}$}
                          &  \bf0.43&\bf(24) & \bf2.08&\bf(37) & \bf13.65&\bf(61) & 244.54&(303) \\ 
 &                        & \bf0.82&\bf(64) & \bf5.89&\bf(123) & \bf46.31&\bf(243) & 523.47&(700) \\
 \cline{2-10}
 &\multirow{2}*{$1.25\times 10^{-3}$} 
                          & \bf0.45&\bf(27) & \bf2.34&\bf(43) & 38.09&(194) & 391.86&(508) \\
 &                        & \bf0.84&\bf(67) & \bf6.17&\bf(128) & 72.27&(393) & 670.94&(907) \\
 \cline{2-10}
 &\multirow{2}*{$2.5\times 10^{-3}$} 
                          & \bf0.51&\bf(32) & 6.54&(132) & 63.74&(326) & 710.77 & (816) \\
 &                        & \bf0.85&\bf(69) & 11.04&(232) & 101.05&(525) & 1160.61 & (1216) \\
 \cline{2-10}
 &\multirow{2}*{$5\times 10^{-3}$} 
                          & 1.26&(96) & 11.15&(223) & 122.05 & (528) & 1863.65 & (2180) \\
 &                        & 1.66&(145) & 16.31&(322) & 188.43 &(726) & 2424.73 & (2571) \\
 \cline{2-10}
 &\multirow{2}*{$10^{-2}$} 
                          &  2.20&(163) & 18.83&(366) & 308.05 & (1459) & 5091.73 & (5724) \\
 &                        & 2.44&(212) & 22.98&(463)  & 386.64 & (1649) & 5979.71 & (6081) \\
 \cline{2-10}
 \cline{1-10}
 \end{tabular}
\end{table}

\begin{figure}[htp]
\begin{center}
\begin{tabular}{cc}
 \includegraphics[width=0.45\textwidth]{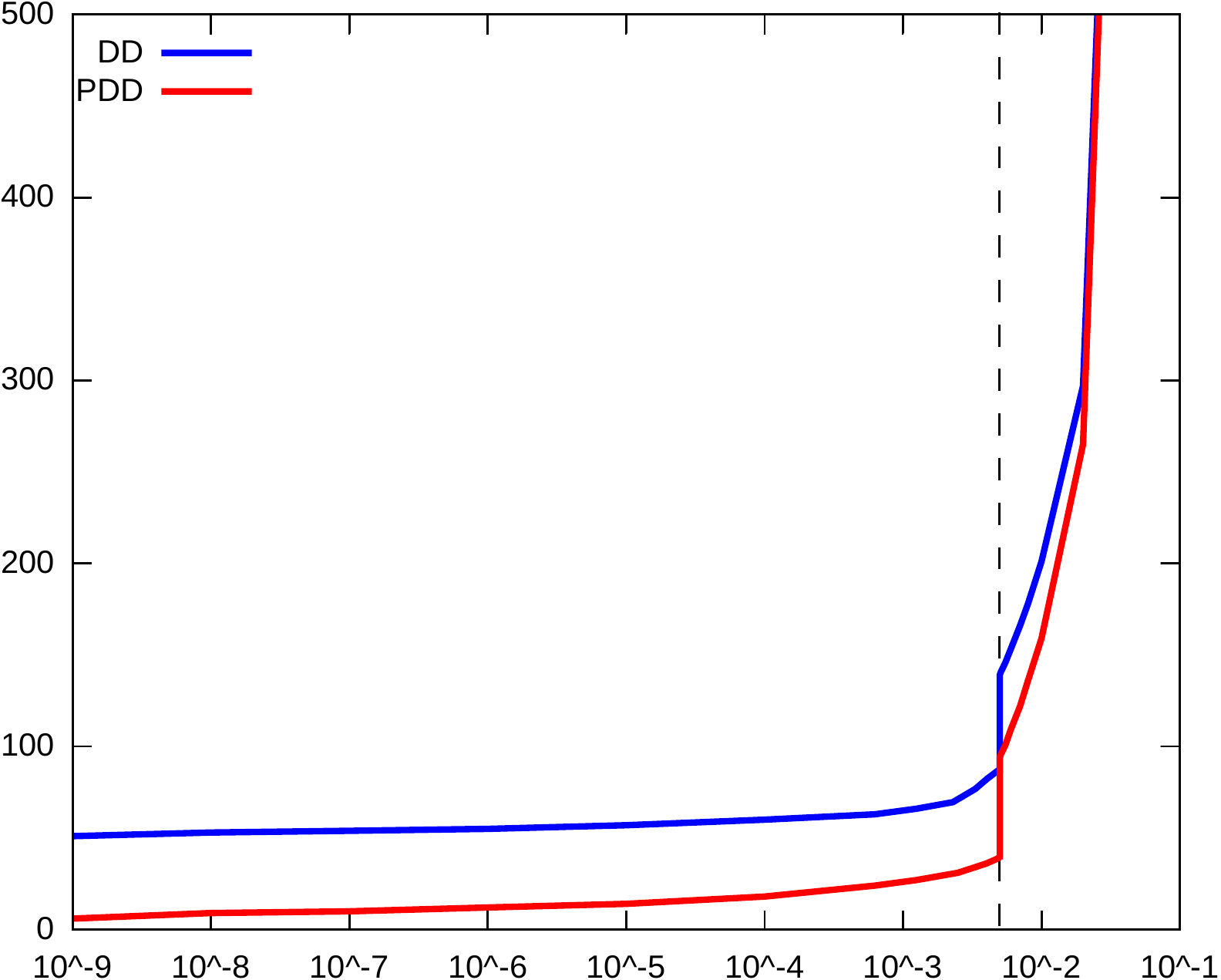} 
 &
 \includegraphics[width=0.45\textwidth]{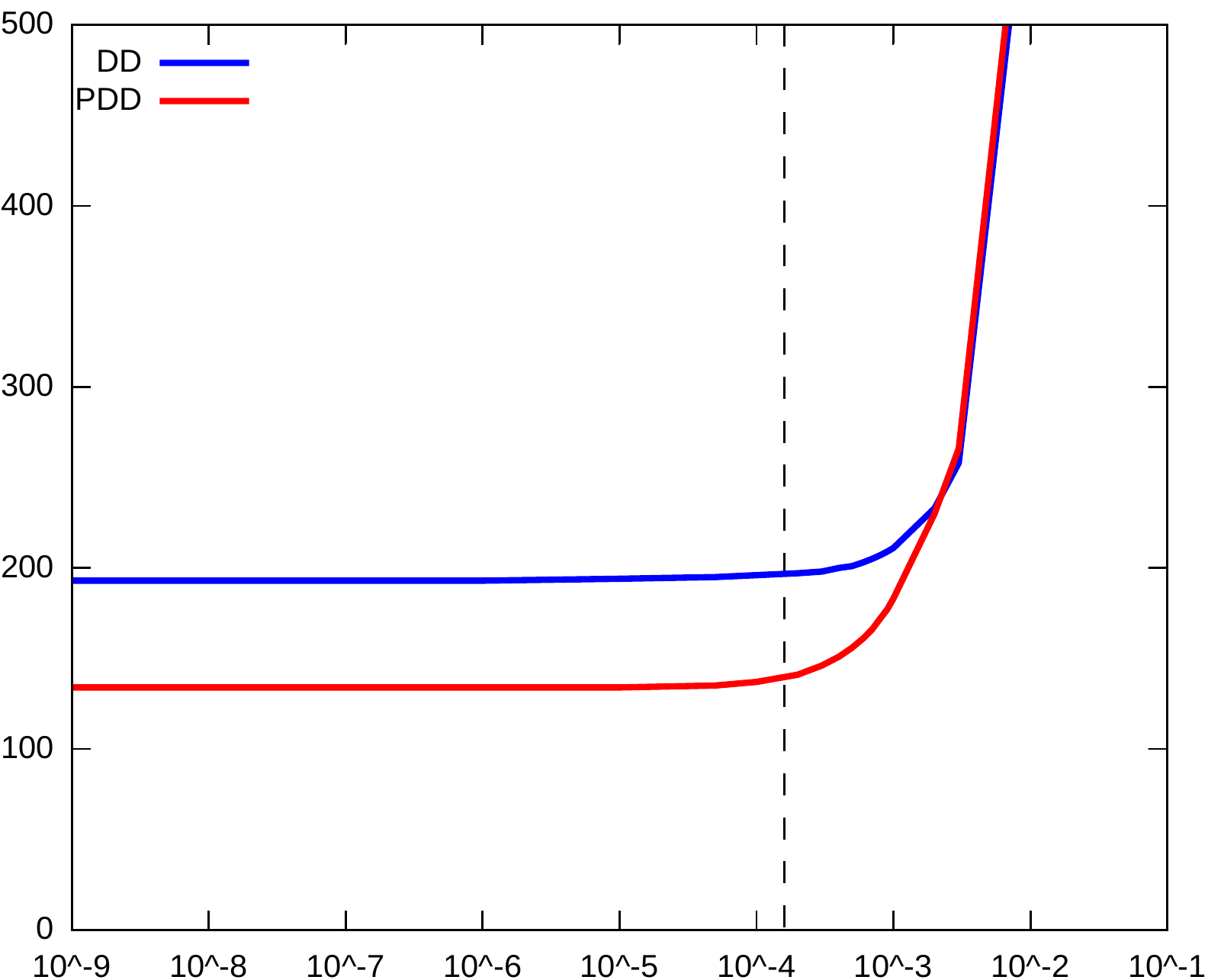} 
 \\
 (a) & (b) 
\end{tabular}
 \caption{
 PDD vs DD in terms of number of iterations (as functions of $\varepsilon$ in 
 logarithmic scale) for $\Delta x=0.02$. 
 Vertical dashed lines identify 
 the thresholds given by the upwind diffusion ball condition: 
 (a) the diffusion-eikonal equation, (b) the Zermelo navigation problem.}\label{eik-graph} 
\end{center}
\end{figure}

Table \ref{table1} refers to the eikonal-diffusion equation. 
The effect of the upwind diffusion ball condition \eqref{hyperbolicity} is completely clear and somehow impressive. 
Indeed, in this test $f_{\min}=f_{\max}=1$ and $\Upsilon=1$, then it reads
$$
\frac{2\varepsilon}{f_{\min}}=\frac{1}{\omega}<\frac{\Delta x}{1+\Upsilon}\qquad\Longleftrightarrow\qquad\varepsilon<\frac{\Delta x}{4}\,.
$$
For the considered fine grids with $\Delta x= 0.02, 0.01, 0.005, 0.0025$ this gives the thresholds 
$\tau_{\Delta x}=5\times 10^{-3}, 2.5\times 10^{-3}, 1.25\times 10^{-3}, 6.25\times 10^{-4}$ respectively. 
We see that both the CPU time and the number of iterations jump 
when the diffusion coefficient 
reaches the corresponding threshold. 
In Figure \ref{eik-graph}a we also show the number of iterations as a 
function of $\varepsilon$ for the case $\Delta x=0.02$ 
with a large number of samples.
This dramatic change of behavior 
is due to the fact that the dynamics is isotropic 
(its module does not depend on the control) 
and constant in speed ($c(x)\equiv 1$). It turns out that, 
in this special case, the upwind diffusion ball condition 
is not only a global but also a local condition, 
the same at each point of the domain. Then, 
for $\varepsilon= 5\times 10^{-3}$ this condition fails 
simultaneously everywhere, producing the observed jump. 
Note that also the standard DD method exhibits the same behavior across the threshold, 
but the performance of the PDD method is much better, 
and this depends, as expected, on the reordering of the nodes in the individual patches. 
 
In Table \ref{table2} we report the results for the Zermelo navigation problem. 
It is easy to see that $f_{\min}=1/6$ and $f_{\max}=3/2$, 
so that $\Upsilon=9$ and the upwind diffusion ball condition reads 
$$
\varepsilon<\frac{\Delta x}{120}\,.
$$
For $\Delta x= 0.02, 0.01, 0.005, 0.0025$ we get the thresholds 
$\tau_{\Delta x}=1.\overline{6}\times 10^{-4}, 8.\overline{3}\times 10^{-5}, 
4.1\overline{6}\times 10^{-5}, 2.08\overline{3}\times 10^{-5}$ respectively. 
We still observe a raising in the number of iterations 
above the corresponding threshold, 
but in this case the effect is moderate, no jump is present (see Figure \ref{eik-graph}b). 
This depends on the fact that the upwind diffusion ball condition 
is here only a global condition, and the threshold estimates 
the worst case scenario of the first points where it fails. 
Moreover, this threshold is not sharp as before. 
Indeed, the value $f_{\min}=1/6$ is attained for instance at 
$(x_1,x_2)=(1,1)$ for $a=e^{i\frac{5\pi}{4}}$, whereas the value $f_{\max}=3/2$ is attained 
for $a=R_\theta \frac{x}{|x|}$ in the limit $x\to 0$. 
Locally the optimal dynamics can do better (and it does!), 
so that, over the threshold, 
we can still have a good performance of the PDD method, 
compared to the standard DD method. 
\begin{table}[htp]
 \centering
 \footnotesize
  \caption{The Zermelo navigation problem: PDD vs DD in terms of CPU time in seconds (and number of iterations) 
  to reach convergence. Values in the first row of each cell concern PDD, the others concern DD. 
  Values in bold are obtained in the regime given by the upwind diffusion ball condition.}\label{table2} 
 \begin{tabular}{|c|c|cc|cc|cc|cc|c|}
 \cline{3-11}
 \multicolumn{2}{c|}{} & \multicolumn{2}{c|}{$100^2$} & \multicolumn{2}{c|}{$200^2$} & \multicolumn{2}{c|}{$400^2$} & \multicolumn{2}{c|}{$800^2$} & $N$\\ \cline{3-11}
 \multicolumn{2}{c|}{} &  \multicolumn{2}{c|}{0.02} & \multicolumn{2}{c|}{0.01} & \multicolumn{2}{c|}{0.005} & \multicolumn{2}{c|}{0.0025} & $\Delta x$ \\ \hline
 \multirow{18}*{$\varepsilon$}
 &\multirow{2}*{$10^{-9}$} & \bf2.94&\bf(132) & \bf16.20&\bf(208) &  \bf96.68&\bf(315)  &  \bf642.71&\bf(522)  \\ 
 &                         & \bf3.24&\bf(191)  & \bf22.52&\bf(325) &  \bf157.74&\bf(577) & \bf1146.91&\bf(1045) \\ 
 \cline{2-10}
 &\multirow{2}*{$10^{-8}$} &  \bf2.97&\bf(132) & \bf16.21&\bf(208) &  \bf96.70&\bf(315) & \bf643.01&\bf(522)\\
 &                         & \bf3.25&\bf(191)  & \bf22.53&\bf(325) & \bf157.75&\bf(577) & \bf1147.27&\bf(1045)\\ 
 \cline{2-10}
 &\multirow{2}*{$10^{-7}$} & \bf2.98&\bf(132) & \bf16.25&\bf(208) & \bf97.14&\bf(315) &  \bf 645.27 & \bf(522) \\ 
 &                        & \bf3.25&\bf(191) & \bf22.54&\bf(325) &  \bf158.15&\bf(577) & \bf 1148.35 & \bf(1046)\\
 \cline{2-10}
 &\multirow{2}*{$10^{-6}$} & \bf2.98&\bf(132) & \bf16.33&\bf(209) & \bf97.23&\bf(316) & \bf 648.84 & \bf (523) \\
 &                        & \bf3.26&\bf(191) & \bf22.57&\bf(325) & \bf158.36&\bf(577) & \bf 1207.03 & \bf (1048) \\
 \cline{2-10}
 &\multirow{2}*{$10^{-5}$} & \bf2.99&\bf(132) & \bf16.30&\bf(209) & \bf 98.91 & \bf (318) & \bf 678.33 & \bf (533) \\
 &                        & \bf 3.26&\bf(191) & \bf22.60&\bf(326) & \bf 159.97 & \bf (581) & \bf 1239.07 & \bf (1061) \\
 \cline{2-10}
 &\multirow{2}*{$2.08\overline3\times 10^{-5}$}
                           & \bf3.01&\bf(132)  & \bf 16.41&\bf(210) & \bf 98.95&\bf (320) & 694.57 & (543)  \\
 &                         & \bf3.27&\bf(191) & \bf 22.67&\bf(327) & \bf 161.16&\bf (583) & 1264.94 & (1075)\\
 \cline{2-10}
 &\multirow{2}*{$4.1\overline6\times 10^{-5}$} 
                          & \bf3.02&\bf(133) & \bf16.76&\bf(212)  & 101.51 & (326)  & 713.82 & (571) \\
 &                        & \bf3.30&\bf(191) & \bf 23.28&\bf(329) & 172.55 & (591) & 1337.13 & (1099) \\
 \cline{2-10}
 &\multirow{2}*{$8.\overline3\times 10^{-5}$} 
                          & \bf3.05&\bf(134) & 16.95&(216)  & 106.20 & (344) & 814.97 & (618) \\
 &                        & \bf3.40&\bf(192) & 23.33&(333) & 167.98 & (604) & 1405.86 & (1149) \\
 \cline{2-10}
  &\multirow{2}*{$1.\overline6\times 10^{-4}$} 
                          & 3.12&(138)       & 17.95&(227)  & 120.27 & (389) & 948.74 & (749) \\
 &                        & 3.35&(194)       & 25.45&(341) & 184.91 & (636) & 1494.67 & (1287) \\
 \cline{2-10}
 &\multirow{2}*{$10^{-3}$} 
                          & 4.08&(187)       & 27.24&(358)  & 242.76 & (774) & 3470.90 & (2827) \\
 &                        & 4.18&(232)       & 33.07&(467) & 312.96 & (996) & 4214.57 & (3261) \\
 \cline{2-10}
 &\multirow{2}*{$5\times 10^{-3}$} 
                          & 8.38&(418) & 111.50 & (1461) & 981.54 & (3221) & 14634.04 & (12105) \\
 &                        & 8.76&(444) & 113.27 & (1507) & 1093.78 & (3360) & 15286.98 & (12348) \\
 \cline{2-10}
 \cline{1-10}
 \end{tabular}
\end{table}

It is important to remark that this is a case where the 
fast-marching method fails, 
due to the strong anisotropy of the problem. 
It follows that the reordering of the nodes according to the 
increasing values of the coarse solution is only 
sub-optimal. Despite we still have an improvement with respect 
to the standard DD method which does not exploit any causality, 
this explains the different performance of the PDD method compared 
to the previous test. 
In particular we see that the number of iterations 
is much larger than that of the eikonal-diffusion equation and 
it remains unchanged for several orders of magnitude of 
$\varepsilon$ (see again Figure \ref{eik-graph}). 
This means that, even for $\varepsilon=0$, 
the dependency between the nodes is already involved 
and very far from the true causality expected in the continuous case for $\Delta x\to 0$.

Finally, we observe in both experiments that, 
as the diffusion coefficient increases, 
the PDD method becomes comparable with the DD method. 
Indeed, the diffusion starts mixing information in the whole domain,
so that the causality is completely lost and the reordering 
turns out to be not only sub-optimal but also useless. 
In this situation the pre-computation step just becomes a 
time wasting procedure.
\section{Conclusions}
We proposed a parallel algorithm for the numerical solution of a class of second order semi-linear equations 
coming from 
stochastic optimal control problems. The new method is based on a dynamic domain decomposition technique, 
and it is a non trivial extension of the patchy domain decomposition method, presented by the authors in \cite{CCFP12} 
for first order Hamilton-Jacobi-Bellman equations. 
We presented a modified semi-Lagrangian local solver, 
in which we removed the self-dependency of the grid nodes induced by the interpolation. This feature, 
combined with fast-marching-like techniques, allowed to accelerate the convergence of the method.
Moreover, we introduced the upwind diffusion ball condition, 
a geometric property that involves the discretization parameters and the data of the problem. 
It guarantees, despite the presence of a diffusion process, 
a qualitative behavior of hyperbolic type for the discretized equation. This translated into an additional 
acceleration for the proposed method.
Several numerical experiments confirmed the ability of the semi-Lagrangian scheme to deal with 
very degenerate diffusion terms and quite general nonlinar equations. 
Finally, we perfomed a comparison with a standard domain decomposition method. We showed that, 
under the regime of the upwind diffusion ball condition and for a sufficiently fine 
grid, the speedup of the new patchy domain decomposition method is remarkable, despite the pre-computation step and the transmission 
conditions. This is the main achievement of this work.

\end{document}